\documentclass[11pt]{article}

\usepackage[letterpaper,margin=1in]{geometry}
\usepackage[utf8]{inputenc}
\usepackage[english]{babel}
\usepackage{enumitem}
\usepackage{ifthen}
\usepackage{subcaption}

\usepackage[colorlinks=true,citecolor=cyan,linkcolor=magenta]{hyperref}
\usepackage{amsmath}
\usepackage{amsfonts}
\usepackage{amssymb}
\usepackage{amsthm}
\usepackage{thmtools}
\usepackage[capitalize,nameinlink]{cleveref} % the hyperlink is only on the number by default. [nameinlink] option includes "Theorem", "Figure", etc in the hyperlink
\crefname{equation}{Equation}{Equations}
\crefname{figure}{Figure}{Figures}
\newtheorem{theorem}{Theorem}[section]
\newtheorem{proposition}[theorem]{Proposition}
\newtheorem{lemma}[theorem]{Lemma}
\newtheorem{corollary}[theorem]{Corollary}

\usepackage{float}
\usepackage{tikz}
\usetikzlibrary{math}
\usetikzlibrary{calc}
\usetikzlibrary{hobby}
\usetikzlibrary{positioning}
\usetikzlibrary{decorations.markings}
\usetikzlibrary{decorations.pathmorphing}
\tikzstyle{vtx}=[circle, draw, fill=black, inner sep=0pt, minimum width=5pt]
\tikzstyle{vtx-white}=[circle, draw, fill=white, inner sep=0pt, minimum width=5pt]
\tikzset{thick/.style={line width=1pt}} % changes 'thick' to set width everywhere
\tikzstyle{matching}=[decorate, decoration={snake, amplitude=3pt, segment length=8pt},line width=1pt]

\newenvironment{restate}[1]
	{\def\thetheorem{#1}\begin{theorem}}
	{\end{theorem}\addtocounter{theorem}{-1}}

\usepackage[colorinlistoftodos,prependcaption,textsize=footnotesize]{todonotes}
\usepackage{xargs}
\newcounter{todocounter}
\newcounter{citation}
\newcounter{prefinaltasks}
\newcommandx{\unsure}[2][1=]{\todo[linecolor=red,backgroundcolor=red!25,bordercolor=red,#1]{#2}}
\newcommandx{\note}[2][1=]{\stepcounter{todocounter}\todo[linecolor=blue,backgroundcolor=blue!25,bordercolor=blue,#1]{\thetodocounter. #2}}
\newcommandx{\fixcitation}[2][1=]{\stepcounter{citation}\todo[caption={\thecitation. fix citation / theorem number},linecolor=green!75!black,backgroundcolor=green!75!black!25,bordercolor=green!75!black,#1]{#2}}
\newcommandx{\nextinline}[2][1=]{\stepcounter{prefinaltasks}\todo[linecolor=violet,backgroundcolor=violet!25,bordercolor=violet,#1]{\theprefinaltasks. #2}}
\newcommandx{\readuntilhere}[2][1=]{\todo[linecolor=green,backgroundcolor=green!25,bordercolor=green]{#2 we've read until here}}

\newcommand{\fig}[3]{
	\begin{figure}
		\centering
		#3
		\caption{#1}
		\label{#2}
	\end{figure}
}
\usepackage{xparse}
\ExplSyntaxOn
\NewDocumentCommand{\subfig}{m m m m}{
	\seq_set_split:Nnn \l_width_seq { , } { #1 }
	\seq_set_split:Nnn \l_caption_seq { , } { #2 }
	\seq_set_split:Nnn \l_label_seq { , } { #3 }
	\seq_set_split:Nnn \l_tikzfile_seq { , } { #4 }

	\int_step_inline:nn { \seq_count:N \l_width_seq }{
		\begin{subfigure}{\seq_item:Nn \l_width_seq { ##1 }\textwidth}
			\centering
			\input{\seq_item:Nn \l_tikzfile_seq { ##1 }}
			\caption{\seq_item:Nn \l_caption_seq { ##1 }}
			\label{\seq_item:Nn \l_label_seq { ##1 }}
		\end{subfigure}%
	}
}
\ExplSyntaxOff

\newcommand{\Xbar}{\overline{X}}
\newcommand{\xbar}{\overline{x}}
\newcommand{\cupdot}{\mathbin{\mathaccent\cdot\cup}}

\makeatletter
\def\moverlay{\mathpalette\mov@rlay}
\def\mov@rlay#1#2{\leavevmode\vtop{%
	\baselineskip\z@skip \lineskiplimit-\maxdimen
	\ialign{\hfil$\m@th#1##$\hfil\cr#2\crcr}}}
\newcommand{\charfusion}[3][\mathord]{
		#1{\ifx#1\mathop\vphantom{#2}\fi
				\mathpalette\mov@rlay{#2\cr#3}
			}
		\ifx#1\mathop\expandafter\displaylimits\fi}
\makeatother

\newcommand{\con}[1]{%
	\def\temp{#1}%
	\def\zero{0}%
	\ifx\temp\empty
		connected%
	\else
		\ifx\temp\zero
			disconnected%
		\else
			\mbox{$#1$-connected}%
		\fi
	\fi
}

\newcommand{\cc}[1]{\mbox{$#1$-connected} cubic graph}

\newcommand{\mc}{matching covered}
\newcommand{\mcg}{\mc~graph}

\newcommand{\sminor}[1]{\mbox{$S$-minor#1}}

\newcommand{\tf}{\mbox{$\theta$-free}}
\newcommand{\tb}{\mbox{$\theta$-based}}
\newcommand{\cfour}{\mathbb{C}_4}
\newtheorem{problem}{Problem}
\crefname{problem}{Problem}{Problems}
\newtheorem{decision-problem}{Decision Problem}
\crefname{decision-problem}{Decision Problem}{Decision Problems}
\newtheorem{statement}{}[theorem]
\crefname{statement}{}{}

\newcommand{\NP}{$\mathsf{NP}$}
\newcommand{\PP}{$\mathsf{P}$}
\newcommand{\coNP}{\mbox{$\mathsf{co}$-\NP}}

\newcommand{\petersen}{\mathbb{P}}
\newcommand{\tftwelve}{K_{3,3}^\Delta}

\newcommand{\yes}{$\sf Yes$}
\newcommand{\no}{$\sf No$}

\title{\tf{} matching covered graphs: characterization and consequences\footnote{Supported by IC\&SR (\url{https://icandsr.iitm.ac.in/}) as well as aRtCS (\url{https://artcs.iitm.ac.in/}).}}
\author{
Rohinee Joshi\footnote{Email: rohi1729@gmail.com}\\
Ajit A. Diwan\footnote{Included as an author posthumously.}\\
\small IIT Bombay
\and
Santhosh Raghul\footnote{Email: cs22d406@cse.iitm.ac.in; Funding source: Ministry of Education (MoE), Government of India.}\\
Nishad Kothari\footnote{Email: nishadkothari@gmail.com}\\
\small IIT Madras
}

\date{%
	\today\\%[1ex]
	Dedicated to beloved Professor Murty
}

\begin{document}

	\maketitle
	\thispagestyle{empty}

	\vspace{-10mm}
	\begin{abstract}
		A nonempty connected graph is {\em matching covered} if each edge lies in some perfect matching; for instance, \mbox{$2$-connected} $3$-regular (aka {\em cubic}) graphs have this property.
One may refer to Lucchesi and Murty's monograph [Perfect Matchings: A Theory of Matching Covered Graphs, 2024] for the extensive literature.
The Ear Decomposition Theorem of Lov\'asz and Plummer [Matching Theory, 1986] implies that every matching covered graph, except $K_2$ and cycles, contains at least one of $\theta$~and~$K_4$ as a `conformal minor' (defined below), where~$\theta$ denotes the loopless cubic graph on two vertices.
Lov\'asz [{\em Combinatorica}, 1983] proved the refinement that every nonbipartite matching covered graph contains at least one of $K_4$ and the triangular prism~$\overline{C_6}$ as a `conformal minor'.
These immediately lead to three problems: characterize `$\theta$-free' (matching covered) graphs, and likewise, `$K_4$-free' and `$\overline{C_6}$-free' ones.

A {\em bisubdivision} of a graph $J$ is any subdivision wherein each edge is replaced by a path of odd length.
% ; any bisubdivision of a matching covered graph (except $K_2$) is also matching covered.
A subgraph $H$ of a graph $G$ is {\em conformal} if $G-V(H)$ has a perfect matching.
A matching covered graph $J$ is a {\em conformal minor} of a matching covered graph $G$ if some conformal subgraph $H$ of $G$ is a bisubdivision of $J$;
otherwise, we say that $G$ is {\em $J$-free}.
For instance, $\theta$ is a conformal minor of $\overline{C_6}$, whereas the Petersen graph is $\theta$-free.
Kothari and Murty [{\em J. Graph Theory}, 2016] used the tight cut decomposition theory to characterize planar graphs that are $K_4$-free, and those that are $\overline{C_6}$-free; the nonplanar case of each problem is open.
In contrast to their work, we exploit a seminal result of Edmonds, Lov\'asz and Pulleyblank [{\em Combinatorica}, 1982], which guarantees the existence of special types of `tight cuts', to obtain a structural $\mathsf{NP}$ characterization of $\theta$-free graphs that immediately places the corresponding decision problem in~$\mathsf{P}$.
The Petersen graph and $K_4$ play key roles in our result. We mention a few consequences.

We deduce that every $\theta$-free graph has at most $2n-2$ edges, where $n$ is its order, and we characterize the tight examples.
Despite being sparse, these graphs are not necessarily planar.
% We provide a characterization of Pfaffian $\theta$-free graphs in the style of Little's [{\em J.~Combin.~Theory Ser.~B}, 1975] co-$\mathsf{NP}$ characterization of Pfaffian bipartite graphs.
In the style of Little's Theorem [{\em J.~Combin.~Theory Ser.~B}, 1975] for bipartite graphs, we characterize Pfaffian $\theta$-free graphs in terms of their forbidden conformal minors.
Using the works of Robertson, Seymour and Thomas [{\em Ann.~of Math.}, 1999], and of McCuaig [{\em The Electronic J.~of Combin.}, 2004], we deduce that the Pfaffian recognition problem is in $\mathsf{P}$ for $\theta$-free graphs (as in the case of bipartite graphs).

It is well-known that deciding whether a cubic graph is $3$-edge-colorable is $\mathsf{NP}$-complete;
for $\theta$-free ones, we provide a characterization of those that are $3$-edge-colorable, and deduce that the corresponding decision problem lies in $\mathsf{P}$.
Note that a conformal cycle is of length either $0$ or $2 \pmod{4}$.
McCuaig [{\em J.~Graph Theory}, 2000] characterized \mbox{$3$-connected} bipartite cubic graphs each of whose conformal cycles is of length $2 \pmod{4}$; the $2$-connected case is open.
We stumbled upon the serendipitous corollary of our main result that each conformal cycle of a $2$-connected cubic graph is of length $0 \pmod{4}$ if and only if it is $\theta$-free.

	\end{abstract}

	% \tableofcontents

	\section{Introduction and summary of results}
		\label{sec:introduction}
		All graphs considered in this paper are loopless; however, we allow multiple/parallel edges.
For general graph theoretic terminology, we follow Bondy and Murty \cite{bomu08}, and for terminology specific to matching theory, we follow Lucchesi and Murty \cite{lumu24}.
This paper may be regarded as a sequel to Kothari and Murty \cite{km16}; however, it is self-contained.
For a graph, we use $n$ for its order and $m$ for its size. 

A graph is {\it matchable} if it has a perfect matching.
Tutte established the following characterization of matchable graphs, where $\mathsf{c_{odd}}(H)$ denotes the number of odd components of a graph~$H$.
By an \emph{odd (even) component}, we mean a component of odd (even) order.

\begin{theorem}{\sc[Tutte's $1$-factor Theorem]}\label{tuttes theorem}
	\newline
	A graph $G$ is matchable if and only if $\mathsf{c_{odd}}(G-S)\leqslant |S|$ for each subset $S$ of $V(G)$. 
\end{theorem}

For a matchable graph~$G$, a set $B \subseteq V(G)$ is a \emph{barrier} if $\mathsf{c_{odd}}(G-B) = |B|$.
A barrier is {\it trivial} if it comprises at most one vertex; otherwise, it is {\it nontrivial}.
Using Tutte's Theorem, one may easily deduce the following.

\begin{proposition}\label{when is G-u-v not matchable}
	For vertices $u$~and~$v$ of a matchable graph~$G$, the graph $G-u-v$ is not matchable if and only if there exists a barrier containing $u$~and~$v$. \qed 
\end{proposition}

A connected nontrivial graph is {\it matching covered} if each of its edges participates in some perfect matching.
There is extensive literature on matching covered graphs; see \cite{lumu24}.
The above proposition immediately yields the characterization of \mc{} graphs stated below. 

\begin{corollary}\label{matchable iff each barrier is stable}
{\sc[Characterization of Matching Covered Graphs]}
\newline
A connected matchable graph is matching covered if and only if each barrier is stable; consequently, if~$B$ is a barrier of a matching covered graph~$G$, then each component of $G-B$ is odd.  
\qed
\end{corollary}

The following fundamental result is a culmination of the works of Kotzig and Lov\'asz; see \cite[Chapter 3]{ lumu24}.

\begin{theorem}{\sc[The Canonical Partition Theorem]}\label{canonical partition theorem}
\newline
	The maximal barriers of a matching covered graph partition its vertex set. 
\end{theorem}

Sch\"onberger (1927) proved that a $3$-regular (aka \emph{cubic}) graph is \mc{} if and only if it is $2$-connected;
this may also be deduced using \cref{tuttes theorem,matchable iff each barrier is stable}.
\cref{fig:smallest 2-conn cubic graphs} shows some such graphs of small order along with their corresponding notations.
Most of these graphs play an important role in the theory of \mcg{}s, including our problem of interest, as we will see in \cref{subsec:tf problem}.

\fig{ A few small \cc2s}{fig:smallest 2-conn cubic graphs}{
	\subfig
	{ 0.15 , 0.175 , 0.15 , 0.225 , 0.225 } % subfigure widths
	{ $\theta$ , $\cfour$ , $K_4$ , $K_{3,3}$ , $\overline{C_6}$ } % subfigure captions
	{fig:theta,fig:cubicc4,fig:k4,fig:k33,fig:c6bar} % subfigure labels
	{ figures/theta.tikz , figures/cubicc4.tikz , figures/k4.tikz , figures/k33.tikz , figures/c6bar.tikz } % tikz file paths
}

Two key aspects of the theory of \mcg{}s are the ear decomposition theory \cite[Chapter 11]{lumu24} and the tight cut decomposition theory \cite[Chapter 4]{lumu24}, each of which is integral to our work.
In particular, the main problem we solve arises from the ear decomposition theory, whereas its solution relies heavily on the tight cut decomposition theory.

		\subsection{\tf{ness} and related problems}
			\label{subsec:tf problem}
			{\it Bisubdividing an edge} refers to the operation of subdividing an edge by inserting an even number of subdivision vertices.
A graph $H$ is a {\it bisubdivision} of a graph $J$ if it may be obtained from $J$ by bisubdividing the members of any (possibly empty) subset of its edge set.
The following is easily observed.

\begin{proposition}\label{every bisubdivision of a mcg is mc}
	Every bisubdivision of a \mcg{}, distinct from $K_2$, is also \mc. \qed
\end{proposition}

We say that a subgraph $H$ of a \mcg{} $G$ is {\it conformal} if $G-V(H)$ is matchable.
Using the theory of ear decompositions (that is discussed in \cref{subsec:tf near-bipartite intro}), Lov\'asz \cite{lova83} proved the following.

\begin{theorem} {\sc [$K_4$--$\overline{C_6}$ Theorem]}\label{k4-c6bar theorem}
% \newline Every nonbipartite matching covered graph either has a conformal bisubdivision of~$K_4$, or has a conformal bisubdivision of $\overline{C_6}$, or possibly both (as a subgraph).
% \newline
% \newline Every nonbipartite matching covered graph has a conformal subgraph that is either a bisubdivision of~$K_4$, or a bisubdivision of $\overline{C_6}$.
% \newline
\newline Every nonbipartite matching covered graph has a conformal subgraph that is a bisubdivision of~$K_4$ or of the triangular prism $\overline{C_6}$.
\end{theorem}

The above theorem inspires the following definition.
A matching covered graph $J$ is a {\it conformal minor} of a matching covered graph~$G$ if the latter has a conformal subgraph~$H$ that is a bisubdivision of~$J$; in this case, for the sake of brevity, we say that~$G$ is {\it $J$-based}; otherwise, $G$ is {\it $J$-free}.
For instance, the Petersen graph $\mathbb{P}$ is \mbox{$K_4$-based} but \mbox{$\overline{C_6}$-free}.
On the other hand, the reader may observe that any plane \mcg{} with two odd faces is \mbox{$K_4$-free}, whence \mbox{$\overline{C_6}$-based}.
It is worth noting that a \mcg{} (for instance, $K_6$) may be \mbox{$K_4$-based} as well as \mbox{$\overline{C_6}$-based}.
The above theorem may be restated as follows: every nonbipartite matching covered graph is either $K_4$-based or $\overline{C_6}$-based, possibly both.
This immediately leads us to the following problems.

\begin{problem}
	\label{problem:K4-free}
	Characterize $K_4$-free matching covered graphs. 
\end{problem}
\begin{problem}
	\label{problem:prism-free}
	Characterize $\overline{C_6}$-free matching covered graphs. 
\end{problem}

Kothari and Murty \cite{koth16,km16} solved the above problems for planar graphs, and their work implies that the corresponding decision problems lie in \PP; however, the nonplanar case of each is open \cite[Problems 12.17 and 12.18]{lumu24}.
We now switch our attention to another result; although reminiscent of \cref{k4-c6bar theorem}, it is an easy consequence of a fundamental result of Lov\'asz and Plummer (\cref{existence of removable ears}).

\begin{theorem} {\sc [$\theta$--$K_4$ Theorem]}\label{theta-k4 theorem}
\newline 
Every matching covered graph, except $K_2$ and cycles, is either \tb{} or $K_4$-based, possibly both.  
\end{theorem}

This once again leads to \cref{problem:K4-free}, but also to the following problem that has surprisingly not received due attention in the past, perhaps because $\theta$ seems too innocent.

\begin{problem}
	\label{problem:theta-free}
	Characterize \tf{} matching covered graphs. 
\end{problem}

For instance, $\overline{C_6}$ is \tb{} whereas the Petersen graph $\petersen$ is \tf{}; see \cref{Petersen is theta free}.
In this paper, we solve \cref{problem:theta-free} by giving a structural description of \tf{} matching covered graphs (\cref{theta-free characterization - inductive version}); the Petersen graph plays a special role.
Our result implies that the corresponding decision problem, stated below, is in \NP.

\begin{decision-problem}\label{dp:tf}
	Given a \mcg{}, decide whether it is \tf{}.
\end{decision-problem}

Additionally, our result leads to a polynomial time algorithm, thus placing the above decision problem in \PP.
Furthermore, using our result, we deduce interesting consequences; for instance, a characterization of Pfaffian \tf{} graphs.

In \cref{subsec:tf graphs characterization}, we shall briefly describe the approach of Kothari and Murty to make progress on \cref{problem:K4-free,problem:prism-free}, and explain why the same does not work for \cref{problem:theta-free}.
Our solution relies heavily on a deep result (\cref{ELP theorem - tight cuts version}) of Edmonds, Lov\'asz and Pulleyblank \cite{elp82}, henceforth abbreviated to ELP.

We conclude by relating \cref{problem:prism-free,problem:theta-free}.
% Note that the triangular prism $\overline{C_6}$ is \tb{}; see \cref{fig:c6bar}.
A spanning (thus, conformal) bisubdivision of $\theta$ in $\overline{C_6}$ is shown in thick pink lines in \cref{fig:c6bar}.
Thus, by transitivity, every \mbox{$\overline{C_6}$-based} \mcg{} is also \tb{}.
However, the converse does not hold;
for instance, wheels (of even order), except $K_4$, are \tb{} but \mbox{$\overline{C_6}$-free}.
Consequently, our solution to \cref{problem:theta-free} may be viewed as characterizing a proper subset of \mbox{$\overline{C_6}$-free} graphs;
thus, solving a special case of \cref{problem:prism-free}.
Henceforth, for the sake of brevity, by \emph{\tf{} (\tb{}) graphs}, we mean \tf{} (\tb{}) \mcg{}s.

		\subsection{Our characterization of \tf{} graphs}
			\label{subsec:tf graphs characterization}
			We start by making a couple of easy observations that we shall find useful later.
The following is easily proved by considering the symmetric difference of two perfect matchings.

\begin{proposition}\label{any two adjacent edges participate in a conformal cycle and C2 is the only nonsimple theta-free mcg}
	Any two adjacent edges of a matching covered graph participate in a conformal cycle.
	Consequently, the cycle $C_2$ is the only \tf{} graph that is not simple.\qed
\end{proposition}

We now state an innocent observation that distinguishes \cref{problem:theta-free} from \cref{problem:K4-free,problem:prism-free}.

\begin{proposition}\label{for simple J: G is J-free iff the underlying simple graph if J-free}
	A matching covered graph $G$ is $J$-free, where $J$ is any simple matching covered graph, if and only if the underlying simple graph of $G$ is also $J$-free. \qed
\end{proposition}

Here, the assumption that $J$ is simple is crucial.
For instance, the graphs shown in \cref{fig:theta,fig:cubicc4} are \tb{}, but their underlying simple graphs are clearly \tf.
% We now continue our discussion on \tf ness.
The proof of the following is implicitly contained in the proof of \cite[Theorem 4.1.6]{lopl86}. 

\begin{proposition}\label{any claw of a bip mcg is part of a conformal bisubdivision of theta}
In a bipartite matching covered graph, any three edges that have a common end participate in a conformal bisubdivision of $\theta$.
\end{proposition}

The above immediately implies the following.

\begin{corollary} {\sc [Bipartite \tf{} Graphs]} \label{the only bipartite theta-free mcgs are k2 and even cycles}
\newline The only bipartite \tf{} graphs are $K_2$ and (even) cycles. 
\end{corollary}

As we shall see, the nonbipartite case of \cref{problem:theta-free} turns out to be far more interesting.
To this end, we need some terminology.
For a subset $X$ of the vertices of a graph $G$, we denote by $\partial(X)$ the {\it cut} associated with~$X$ --- that is, the set of edges that have one end in $X$ and the other end in $\overline{X}:=V(G)-X$.
We refer to $X$ and $\overline{X}$ as the {\it shores} of the cut $C:=\partial(X)$.
For a vertex~$v$, we simplify the notation to $\partial(v):=\partial(\{v\})$; such a cut is called {\it trivial}.
In the same spirit, for a subgraph $H$, we simplify $\partial(V(H))$ to $\partial(H)$.
When $G$ is of even order, a cut $C$ is {\it odd} if both shores have an odd number of vertices; otherwise, $C$ is {\it even}.

For a cut $C:=\partial(X)$ of a graph $G$, the graph obtained by identifying each vertex in $\overline{X}$ to a single vertex $\overline{x}$ is denoted by $G/(\overline{X}\rightarrow \overline{x})$, or simply by $G/\overline{X}$.
We say that $G/\overline{X}$ is obtained from $G$ by \emph{shrinking the shore} $\overline{X}$ of the cut $C$.
The two graphs $G/X$ and $G/\overline{X}$ are called the {\it $C$-contractions} of $G$.
\cref{fig:cuts and contractions} shows an example of a cut $C$, and the corresponding $C$-contractions.

\fig{A cut $C$ and the corresponding $C$-contractions}{fig:cuts and contractions}{
	\subfig{0.25,0.42,0.33}{$G/(X \to x)$,Tight cuts in a \mcg{} $G$,$T_6^+:=G/(\Xbar \to \xbar)$}{fig:tight cut contraction 1,fig:tight cut example,fig:tight cut contraction 2}{figures/tight-cut-contraction-example-1.tikz,figures/tight-cut-example.tikz,figures/tight-cut-contraction-example-2.tikz}
}

Next, we discuss special types of cuts that play a crucial role in the theory of \mcg{}s \cite[Chapter 4]{lumu24}; in particular, they were also employed by Kothari and Murty in \cite{km16} to solve the planar case of \cref{problem:K4-free,problem:prism-free}.

			\subsubsection{Tight cuts and the tight cut decomposition procedure}
				\label{subsubsec:tight cuts}
				A cut $C:=\partial(X)$ of a matching covered graph $G$ is a {\it tight cut} if $|C\cap M|=1$ for each perfect matching $M$; for instance, both cuts $C$ and $D$ shown in \cref{fig:tight cut example} are tight.
It is easily seen that if $C$ is a tight cut, then both $C$-contractions (aka \emph{tight cut contractions}) are matching covered.
In particular, if $C$ is nontrivial, then the $C$-contractions are of smaller order;
naturally, if either of them has a nontrivial tight cut, then one may obtain two smaller graphs in the same manner.
One may perform this process repeatedly until one obtains a list of matching covered graphs --- each of which is free of nontrivial tight cuts.
This process is called the {\it tight cut decomposition procedure}.
For instance, the reader may verify that any application of the tight cut decomposition procedure to the graph shown in \cref{fig:tight cut example} yields three copies of $K_4$ (up to multiple edges).

A matching covered graph that is free of nontrivial tight cuts is called a {\it brace} if it is bipartite, or a {\it brick} if it is nonbipartite.
For instance, the graphs shown in \cref{fig:smallest 2-conn cubic graphs} are all of the bricks and braces of order at most six that are cubic.
% \cref{fig:smallest 2-conn cubic graphs,fig:bipartite ear decomposition,fig:nonbipartite ear decomposition}
% are bricks and braces.
In general, a matching covered graph may admit several applications of the tight cut decomposition procedure.
However, Lov\'asz \cite{lova87} proved the following remarkable result.

\begin{theorem}{\sc [Unique Tight Cut Decomposition Theorem]}\label{unique tcd theorem}
\newline
Any two applications of the tight cut decomposition procedure to a matching covered graph yield the same list of bricks and braces, up to multiplicities of edges.
\end{theorem}

In light of the above, by \emph{bricks and braces of} (a matching covered graph)~$G$, we mean the underlying simple graphs produced by any application of the tight cut decomposition procedure.
In particular, the above result implies that the number of bricks of~$G$, denoted by $b(G)$, is an invariant;
it plays a key role in matching theory; see~\cite{lumu24}.
For instance, $b(G)=3$ for the graph~$G$ shown in \cref{fig:tight cut example}.
It is worth noting that $G$ is bipartite if and only if $b(G)=0$.

A matchable graph~$G$ is {\it bicritical} if $G-u-v$ is matchable for each pair of distinct vertices~$u$~and~$v$. 
We remark that a matching covered graph has zero braces if and only if it is bicritical and has order four or more.
This is one of the reasons as to why most authors, including \cite{km16,lumu24}, require bicritical graphs to have order at least four;
however, we shall find it convenient to not impose this restriction.
As per our definition, $K_2$ (up to multiple edges) is the only bicritical graph that has precisely two vertices, and the only one that is bipartite.
We now proceed to state a characterization of bricks due to Edmonds, Lov\'asz and Pulleyblank~\cite{elp82}.

\begin{theorem}{\sc[Characterization of Bricks]}\label{ELP theorem - bricks version}\newline
	Bricks are precisely the $3$-connected bicritical graphs of order four or more.
\end{theorem}

In their pursuit of solving \cref{problem:K4-free,problem:prism-free}, Kothari and Murty \cite{km16} proved the following result as a first step.

\begin{theorem}\label{for a cubic brick J: G is J-free iff each tight cut contraction is J-free}
	For any cubic brick $J$, a \mcg{} $G$ is $J$-free if and only if both of its $C$-contractions are $J$-free, where $C$ is any tight cut of $G$.
\end{theorem}

Since braces are bipartite and cubic bricks are simple, the above combined with \cref{for simple J: G is J-free iff the underlying simple graph if J-free} implies the following.

\begin{corollary}\label{for a cubic brick J: G is J-free iff each brick of G is J-free}
	For any cubic brick $J$, a matching covered graph~$G$ is $J$-free if and only if each of its bricks is $J$-free. \qed
\end{corollary}

The bricks (and braces) of a matching covered graph may be computed in polynomial time \cite[Exercise 5.4.5]{lumu24}.
Hence, due to the above corollary, it suffices to characterize $K_4$-free and $\overline{C_6}$-free bricks in order to solve \cref{problem:K4-free,problem:prism-free}, respectively.

Naturally, one is tempted to wonder whether the above approach may be employed to solve \cref{problem:theta-free}.
Unfortunately, however, neither \cref{for simple J: G is J-free iff the underlying simple graph if J-free} (as discussed earlier) nor \cref{for a cubic brick J: G is J-free iff each tight cut contraction is J-free} holds for $J=\theta$.
To see the latter, the reader may easily verify that $T_6$, shown in \cref{fig:T6-first}, is \tf{}, whereas each of its $C$-contractions is isomorphic to the \tb{} graph shown in \cref{fig:C-contraction of T6}.

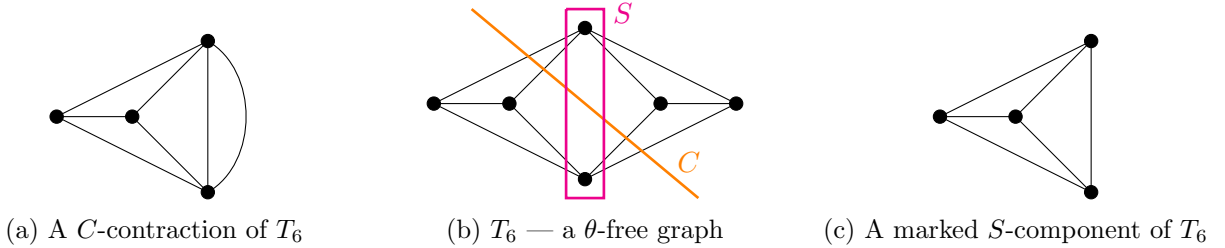
\fig{$C$-contraction and marked $S$-component of $T_6$}{fig:T6ill}{
	\subfig{0.31,0.38,0.31}{A $C$-contraction of $T_6$,$T_6$ --- a \tf{} graph,A marked $S$-component of $T_6$}{fig:C-contraction of T6,fig:T6-first,fig:marked S-component of T6}{figures/C-contraction-of-T6.tikz,figures/T6-2-sep-cut.tikz,figures/marked-S-component-of-T6.tikz}
}

% One of the issues, as we realized, is that \cref{for a cubic brick J: G is J-free iff each tight cut contraction is J-free} deals with arbitrary tight cuts.
% Luckily, a result of Edmonds, Lov\'asz and Pulleyblank, that is equivalent to \cref{ELP theorem - bricks version},
% proves the existence of special types of nontrivial tight cuts called `barrier cuts' and `\mbox{$2$-separation} cuts' (whenever a nontrivial tight cut exists).
% As we shall see soon, their result comes to our rescue in our pursuit of solving \cref{problem:theta-free}. 

On one hand, \cref{for a cubic brick J: G is J-free iff each tight cut contraction is J-free} is quite strong as it deals with arbitrary tight cuts.
On the other hand, it has its limitations as it is applicable only to cubic bricks; in particular, it is not applicable to the graph $\theta$ of our interest.
This realization led us to wonder whether dealing with special types of nontrivial tight cuts would help us in making progress towards solving \cref{problem:theta-free}.
However, for such an approach to work, one would need to show the existence of such tight cuts (whenever a nontrivial tight cut exists).
Luckily, Edmonds, Lov\'asz and Pulleyblank \cite{elp82} have already done the heavy lifting, and their result (equivalent to \cref{ELP theorem - bricks version}) comes to our rescue as discussed next.

			\subsubsection{Leveraging ELP cuts to solve \cref{problem:theta-free}}
				\label{subsubsec:ELP cuts}
				Let $G$ be any \mcg{} and let $B$ denote a barrier of $G$.
It is easily observed that for any (odd) component $L$ of $G-B$, the cut $\partial_G(L)$ is tight; a tight cut that arises from a barrier, in this manner, is a {\it barrier cut}.
For an illustration, see the graphs shown in \cref{fig:barrier}.

\fig{Barrier cuts}{fig:barrier}{
	\subfig{0.28,0.325,0.385}{A \tb{} graph,Another \tb{} graph,$\tftwelve$ --- a \tf{} graph}{fig:barriertb1,fig:barriertb2,fig:barriertf}{figures/k33-splice-k4.tikz,figures/k33-splice-2-k4s.tikz,figures/k33-splice-3-k4s.tikz}
}

The (bipartite matching covered) graph $H$ obtained from $G$, by shrinking each (odd) component of $G-B$ to a single vertex, is the \emph{core of $G$ with respect to $B$}, or simply its \emph{$B$-core}.
If $L$ is any component of $G-B$, the graph obtained by shrinking $\overline{V(L)}:=V(G)-V(L)$ into a single vertex is a \emph{barrier fragment of~$G$ with respect to $B$}, or simply a \mbox{\emph{$B$-fragment}}.
Note that there are exactly $|B|$ many \mbox{$B$-fragments}, and that their edge sets comprise a partition of $E(G)$.
For instance, for the graph shown in \cref{fig:barriertb2}, its $B$-core is $K_{3,3}$ and of its three \emph{$B$-fragments}, two are (isomorphic to) $K_4$ and one is $\theta$.

For the \tf{} graph shown in \cref{fig:barriertf}, one of its $C$-contractions is the \tb{} graph shown in \cref{fig:barriertb2}.
This shows that \cref{for a cubic brick J: G is J-free iff each tight cut contraction is J-free} does not hold for~$J=\theta$ even if $C$ is a barrier cut.
Interestingly, it turns out that instead of considering a particular barrier cut arising from a barrier~$B$, one should pay attention to all $B$-fragments (and ignore the $B$-core), as stated below.

\begin{theorem}{\sc [Recursing along a Barrier]}\label{J-freeness across a barrier where J is bicritical cubic}\newline
	For any bicritical cubic graph $J$, a \mcg{} $G$ is $J$-free if and only if each of its $B$-fragments is $J$-free, where $B$ is any barrier of $G$.
\end{theorem}

The above result is one of our key steps towards solving \cref{problem:theta-free}; a proof appears in \cref{sec:barrier}.
The reader may find it instructive to use it to explain why the first two graphs shown in \cref{fig:barrier} are \tb{} whereas the last one $\tftwelve$ is surprisingly \tf{}; interestingly, $\tftwelve$ has bisubdivisions of $\theta$, none of which are conformal.

It is easily proved that, in a bipartite matching covered graph, every tight cut is a barrier cut; see \cite[Theorem 4.8]{lumu24}.
However, this does not hold for nonbipartite graphs; we now proceed to describe another type of tight cuts.
For a \mcg{} $G$, a~\mbox{{\it $2$-separation}} is any \mbox{$2$-vertex-cut} $S:=\{u,v\}$ that is not a barrier;
consequently, each component of $G-S$ is even.
Now, let $L'$ denote the union of any nonempty proper subset of the (even) components of $G-S$, and observe that the nontrivial cut $\partial_G(V(L') \cup \{u\})$ is tight;
a tight cut that arises from a \mbox{$2$-separation}, in this manner, is called a \mbox{{\it $2$-separation cut}}.
See \cref{fig:2sep} for a few examples.

Let us revisit the \tf{} graph $T_6$ whose both \mbox{$C$-contractions}, isomorphic to the graph shown in \cref{fig:C-contraction of T6}, are \tb{}.
This shows that \cref{for a cubic brick J: G is J-free iff each tight cut contraction is J-free} does not hold for $J=\theta$ even if $C$ is a $2$-separation cut.
Interestingly, it turns out that instead of considering the \mbox{$C$-contractions}, one should consider the marked $S$-components as defined in \cite[Section 9.4]{bomu08} and restated here.
Let $S$ be any $2$-separation of a \mcg{} $G$ and let $L'$ denote any (even) component of $G-S$.
The graph $L$ obtained from the induced subgraph $G[V(L')\cup S]$ by adding a new edge joining~$u$~and~$v$, called its {\it marker edge}, is the {\it marked $S$-component of~$G$ corresponding to $L$}.
For instance, the graph $T_6$ has two marked $S$-components each of which is isomorphic to $K_4$ (see \cref{fig:marked S-component of T6}), whereas the graph shown in \cref{fig:2septb3} has three marked $S$-components.

\fig{$2$-separations}{fig:2sep}{
	\subfig{0.3,0.3,0.3}{$T_6$ --- a \tf{} graph,A \tb{} graph,Another \tb{} graph}{fig:T6,fig:2septb,fig:2septb3}%
	{figures/2-sep-T6.tikz,figures/2-sep-theta-based-example-1.tikz,figures/2-sep-theta-based-example-2.tikz}
}

We invite the reader to compare $L$, defined in the preceding paragraph, with the tight cut contraction $G':=G/(\overline{V(L)\cup \{u\}}\rightarrow v)$ and observe that they are the same, except possibly for the (nonzero) number of edges joining $u$~and~$v$.
As per \cref{for simple J: G is J-free iff the underlying simple graph if J-free}, $L'$ is $J$-free if and only if $G'$ is $J$-free, where $J$ is any simple matching covered graph.
However, the graph of our interest $J = \theta$ is not simple, and considering the marked $S$-components instead of the tight cut contractions turns out to be crucial.
For instance, both marked $S$-components of $T_6$ are $K_4$; more importantly, they are \tf.
The above discussion about $G'$ and $L$ immediately proves statements $(i)$ and $(ii)$ of the result stated below, and statement $(iii)$ follows from \cref{any two adjacent edges participate in a conformal cycle and C2 is the only nonsimple theta-free mcg}.

\begin{proposition}\label{observations on 2-separations}
	Let $S:=\{u,v\}$ denote a $2$-separation of a \mcg{} $G$,
	let $L'$ be any component of $G-S$,
	and let $L$ be the corresponding marked $S$-component.
	Then: (i)~$L$ is matching covered, (ii)~$L'$ is matchable, and (iii)~there exists a conformal (odd) $uv$-path in~$L$ that avoids its marker edge. \qed
\end{proposition}

The following two results comprise our second step towards solving \cref{problem:theta-free}.
The first of these, which follows immediately from \cref{observations on 2-separations} $(ii)$ and $(iii)$, is stated below; \cref{fig:2septb3} illustrates an example.

\begin{proposition}\label{G with a 2-separation S is theta-based if G-S has at least three components}
	Every \mcg{} that has a $2$-separation, whose deletion yields three or more components, is \tb{}. \qed
\end{proposition}

In light of the above, the technical assumption at the end of the next statement is quite innocuous, and in fact necessary.

\begin{theorem}{\sc [Recursing across a $2$-separation]}\label{J-freeness across a 2-separation where J is 3-connected cubic}\newline
	For any $3$-connected cubic graph $J$, a \mcg{} $G$ is $J$-free if and only if each of its each marked $S$-components is $J$-free, where $S$ is any $2$-separation of $G$ --- under the additional assumption that if $J=\theta$ then $G-S$ has precisely two components.
\end{theorem}

A proof appears in \cref{sec:2sep}.
The reader may find it instructive to use the above two results to explain why $T_6$ (shown in \cref{fig:T6}) is \tf;
whereas the other graphs shown in \cref{fig:2sep} and the graph $T_6^+$ (obtained from $T_6$ by adding an edge joining the two vertices of degree four), shown in \cref{fig:tight cut contraction 2}, are all \tb{}.
Having dealt with barrier cuts and $2$-separation cuts via \cref{J-freeness across a barrier where J is bicritical cubic,J-freeness across a 2-separation where J is 3-connected cubic}, respectively, we now state the seminal result of ELP that comes to our rescue as mentioned at the end of \cref{subsubsec:tight cuts}.

\begin{theorem}{\sc[ELP Theorem]}\label{ELP theorem - tight cuts version}\newline
	Every matching covered graph, that has a nontrivial tight cut, has one that is either a barrier cut or a $2$-separation cut, possibly both.
\end{theorem}

In light of the above, by an {\it ELP cut}, we mean a nontrivial tight cut that is either a barrier cut or a $2$-separation cut.
Not all nontrivial tight cuts are ELP cuts;
for instance, the reader may verify that the tight cut~$C$, shown in \cref{fig:tight cut example}, is not an ELP cut, whereas $D$ is.

We now invite the reader to observe that \cref{J-freeness across a barrier where J is bicritical cubic,J-freeness across a 2-separation where J is 3-connected cubic}, along with the ELP Theorem (\ref{ELP theorem - tight cuts version}), imply \cref{for a cubic brick J: G is J-free iff each brick of G is J-free} of Kothari and Murty.
Furthermore, for any bicritical cubic graph~$J$, we infer that it suffices to characterize $J$-free bricks in order to describe all $J$-free matching covered graphs --- in the case of $J=\theta$, one needs to additionally invoke \cref{G with a 2-separation S is theta-based if G-S has at least three components,the only bipartite theta-free mcgs are k2 and even cycles}.
Ergo, to solve \cref{problem:theta-free}, it only remains to characterize the \tf{} bricks, and this brings us to presenting the final piece of the puzzle.
% \unsure{ question: does it also imply \cref{for a cubic brick J: G is J-free iff each tight cut contraction is J-free} \newline answer: yes, but only using \cref{for a cubic brick J: G is J-free iff each brick of G is J-free} and i don't see a way around it}

\begin{restatable}{theorem}{tfBricks}
{\sc[Characterization of \tf{} Bricks]}\label{K4 and Petersen are the only theta-free bricks}\newline
	The only \tf{} bricks are $K_4$ and the Petersen graph $\mathbb{P}$.
\end{restatable}

\fig{The Petersen graph}{fig:petersen}{\begin{tikzpicture}[scale=1.5]

\tikzmath{
	\a = 0.78;
	\b = 1.5;
}

\node[vtx] (v0) at (0,0) {};
\foreach \i in {1,...,9} {
	\node[vtx] (v\i) at (10+\i*40:\b) {};
}

% \foreach \i in {1,...,8} {
% 	\tikzmath{\j=int(\i+1);}
% 	\draw (v\i) -- (v\j);
% }
% \draw (v1) -- (v9);

\draw (0,0) circle (\b);

\draw (v0) -- (v2);
\draw (v0) -- (v5);
\draw (v0) -- (v8);
\draw (v3) to[bend right=30] (v7);
\draw (v4) to[bend left=30] (v9);
\draw (v6) to[bend right=30] (v1);

% \begin{scope}[shift={(4,0)}]
% \node[vtx] (v0) at (0,0) {};
% \foreach \i in {1,...,9} {
% 	\node[vtx] (v\i) at (10+\i*40:\b) {};
% }

% \foreach \i in {1,...,8} {
% 	\tikzmath{\j=int(\i+1);}
% 	\draw (v\i) -- (v\j);
% }
% \draw (v1) -- (v9);

% \draw (v0) -- (v2);
% \draw (v0) -- (v5);
% \draw (v0) -- (v8);
% \draw (v3) to[bend right=22] (v7);
% \draw (v4) to[bend left=22] (v9);
% \draw (v6) to[bend right=22] (v1);
% \end{scope}

% \begin{scope}[shift={(8,0)}]

% \foreach \i in {0,...,4} {
% 	\node[vtx] (v\i) at (18+\i*72:\a) {};
% 	\node[vtx] (u\i) at (18+\i*72:\b) {};
% }

% \draw (v0) -- (u0);
% \draw (v1) -- (u1);
% \draw (v2) -- (u2);
% \draw (v3) -- (u3);
% \draw (v4) -- (u4);

% \draw (u0) -- (u1);
% \draw (u1) -- (u2);
% \draw (u2) -- (u3);
% \draw (u3) -- (u4);
% \draw (u4) -- (u0);

% \draw (v0) -- (v2);
% \draw (v2) -- (v4);
% \draw (v4) -- (v1);
% \draw (v1) -- (v3);
% \draw (v3) -- (v0);

% \end{scope}

% \begin{scope}[shift={(4,-4)}]

% \foreach \i in {0,...,4} {
% 	\node[vtx] (u\i) at (18+\i*72:0.6) {};
% 	\node[vtx] (v\i) at (18+\i*72:1.5) {};
% }

% \draw (v0) -- (u0);
% \draw (v1) -- (u1);
% \draw (v2) -- (u2);
% \draw (v3) -- (u3);
% \draw (v4) -- (u4);

% \draw (u0) to[bend left=30] (u1);
% \draw (u1) to[bend left=30] (u2);
% \draw (u2) to[bend left=30] (u3);
% \draw (u3) to[bend left=30] (u4);
% \draw (u4) to[bend left=30] (u0);

% \draw (v0) to[bend right=24] (v2);
% \draw (v2) to[bend right=24] (v4);
% \draw (v4) to[bend right=24] (v1);
% \draw (v1) to[bend right=24] (v3);
% \draw (v3) to[bend right=24] (v0);

% \end{scope}

\end{tikzpicture}}

We use the Thin Edge Theorem (\ref{CLM thin edge theorem}) of Carvalho, Lucchesi and Murty \cite{clm06} to prove the above in \cref{sec:bricks}.
Putting all of the aforementioned results (\ref{the only bipartite theta-free mcgs are k2 and even cycles}, \ref{J-freeness across a barrier where J is bicritical cubic}, \ref{G with a 2-separation S is theta-based if G-S has at least three components}, \ref{J-freeness across a 2-separation where J is 3-connected cubic}, \ref{ELP theorem - tight cuts version} and \ref{K4 and Petersen are the only theta-free bricks}) together, we now deduce our description of \tf{} graphs.

\begin{theorem}{\sc[Characterization of \tf{} Graphs - Inductive Version]}\label{theta-free characterization - inductive version}\newline
	A matching covered graph~$G$ is \tf{} if and only if at least one of the following holds: 
	\begin{enumerate}[label=(\roman*)]
		\item $G$ is $K_2, C_2, K_4$ or the Petersen graph $\mathbb{P}$, or
		\item $G$ has a nontrivial barrier and for any such barrier, say $B$, each $B$-fragment is \tf,~or
		\item $G$ has a $2$-separation and for any such $2$-separation, say $S$, there are precisely two marked $S$-components and both of them are \tf. \qed
	\end{enumerate}
\end{theorem}

The significance of the ELP Theorem (\ref{ELP theorem - tight cuts version}) to our work, especially to the above result, cannot be overstated.
Ergo, we make some historical remarks pertaining to their work.
ELP proved their result in the pursuit of computing the dimension of the matching space; see \cite[Section 6.3]{lumu24}.
It is worth noting that Naddef \cite{nadd82} computed the same using different techniques.
Surprisingly, the original proof of the ELP Theorem relied on linear programming duality.
Szigeti~\cite{szig02} gave the first graph-theoretical proof;
inspired by his work, Carvalho, Lucchesi and Murty~\cite{clm18} gave yet another proof that relies heavily on the theory developed by them; see \cite[Section 5.2]{lumu24}. They~\cite{clm18} conjectured a stronger version of the ELP Theorem, and this was proved recently by Chen, Feng, Lu, Lucchesi and Zhang \cite{cfllz21}; see \cite[Section 5.5.1]{lumu24}.
We conclude this section with the following easy consequence of \cref{theta-free characterization - inductive version} and Lov\'asz's Unique Tight Cut Decomposition Theorem (\ref{unique tcd theorem}).

\begin{corollary}
	\label{each brick of a tf graph is either K4 or the Petersen graph}
	Each brick of a \tf{} graph is either $K_4$ or the Petersen graph $\petersen$.
	\qed
\end{corollary}

In the next section, we state some interesting corollaries of our main result (\cref{theta-free characterization - inductive version}).

		\subsection{Consequences of our characterization of \tf{} graphs}
			\label{subsec:consequences}
			Our Main Theorem (\ref{theta-free characterization - inductive version}) implies that \cref{dp:tf} lies in \NP;
we also state an equivalent version of the above, from a generation viewpoint, in \cref{theta-free characterization - generation version}.
We now proceed to describe a polynomial time algorithm, thus placing \cref{dp:tf} in \PP;
we shall use the fact that the canonical partition (recall \cref{canonical partition theorem}) of a \mcg{} may be computed in polynomial time using \cref{when is G-u-v not matchable}.

Let $G$ be the given \mcg.
If $G$ is one of $K_2,C_2,K_4$ and the Petersen graph~$\petersen$, return \yes{}.
Next, compute the canonical partition of~$G$.
If any part, say $B$, has cardinality two or more, recurse on each of the $B$-fragments and return \yes{} if each recursive call returns \yes{}; otherwise, return \no{}.
(If each part is trivial then, by \cref{when is G-u-v not matchable}, $G$ is bicritical.)
Next, check whether $G$ has a $2$-vertex-cut $S$; if not, return \no{}.
(Since $G$ is bicritical, $S$ is a $2$-separation.)
Next, compute the marked $S$-components of $G$; if there are three or more, return \no{}.
Finally, recurse on both marked $S$-components and return \yes{} if each recursive call returns \yes{}; otherwise, return \no{}.

In what follows, we discuss other consequences of our main result.
The first of these confirms the intuition that \tf{} graphs ought to be sparse.

			\subsubsection{Sparsity and bicriticality}
				\label{subsubsec:sparsity and bicriticality intro}
				Intuitively, one would expect an upper bound on the size of a \tf{} graph.
This turns out to be true, and an easy consequence of \cref{theta-free characterization - inductive version}.
We start by defining an operation that may be viewed as the inverse of considering the marked $S$-components with respect to a $2$-separation $S$.

% Recall the $K_2$-sum operation, discussed in \cref{subsubsec:sparsity and bicriticality intro}, and note that $S:=\{u,v\}$ is a $2$-separation of $G$ and that $G$ has precisely two marked $S$-components.

Let $G_1$ and $G_2$ be disjoint graphs, and for $i \in \{1,2\}$, let $e_i:=u_iv_i$ be an edge of $G_i$;
the graph obtained from $(G_1-e_1) \cupdot (G_2-e_2)$ by
identifying~$u_1$~and~$u_2$ into a single vertex~$u$, and
likewise $v_1$~and~$v_2$ into a single vertex~$v$,
is called a {\it $K_2$-sum of $G_1$~and~$G_2$
at the edges $e_1$~and~$e_2$}.
For instance, the graph $T_6$, shown in \cref{fig:T6}, is the $K_2$-sum of two copies of $K_4$.
% We remark that this operation is reminiscent of the $C_4$-sum (aka $4$-cycle sum) operation that was crucial in the characterization of \mbox{$K_{3,3}$-free} bipartite \mcg{}s, or equivalently Pfaffian bipartite \mcg{}s, by Robertson, Seymour and Thomas \cite{rst99}, and McCuaig~\cite{mccu04}. 
The reader may easily verify that any $K_2$-sum of \mcg{}s, neither of which is isomorphic to $K_2$, is also \mc{}.
We now define two graph families using the \mbox{$K_2$-sum} operation; one of them, denoted by $\mathcal{T}_0$, is a proper subset of the other, denoted by $\mathcal{T}$.

The brace~$C_2$, the bricks $K_4$ and the Petersen graph $\mathbb{P}$ belong to $\mathcal{T}$,
% and any graph obtained from two members of $\mathcal{T}$ using the \mbox{$K_2$-sum} operation belongs to $\mathcal{T}$.
and any \mbox{$K_2$-sum} of two members of $\mathcal{T}$ also belongs to $\mathcal{T}$.
The brace~$C_2$ and the brick $K_4$ belong to $\mathcal{T}_0$,
% and any graph obtained from two members of~$\mathcal{T}_0$ using the \mbox{$K_2$-sum} operation belongs to $\mathcal{T}_0$.
and any \mbox{$K_2$-sum} of two members of~$\mathcal{T}_0$ also belongs to $\mathcal{T}_0$.
Observe that the inclusion of $C_2$ does not result in any additional elements since the $K_2$-sum of $C_2$ and any graph $G$ is simply $G$.
The graph $T_6$ shown in \cref{fig:T6} is the third smallest member of $\mathcal{T}_0$ as well as of~$\mathcal{T}$.
Recall that $b$ denotes the number of bricks of a \mcg{}.
We prove the following in \cref{subsec:sparsity and bicriticality}.

\begin{restatable}{theorem}{edgeBoundTheorem}
\label{edge bound on theta-free graphs}
	Every \tf{} graph $G$ satisfies:
	\begin{enumerate}[label=(\roman*)]
		\item $m\leqslant \frac{3n}{2}+b-1$, and equality holds if and only if $G\in \mathcal{T}$,
		\item $b\leqslant \frac{n}{2}-1$, and equality holds if and only if $G\in\mathcal{T}_0$, and
		\item consequently, $m\leqslant 2n-2$, and equality holds if and only if $G\in\mathcal{T}_0$.
	\end{enumerate}
\end{restatable}

In fact, it is easy to prove that the upper bound in statement $(ii)$ of the above holds for all \mcg{}s; see \cite[Exercise 4.4.3]{lumu24}.
However, there are \tb{} graphs that also satisfy this bound with equality; one such example is $T_6^+$, shown in \cref{fig:tight cut contraction 2}.

In light of the upper bound in statement $(iii)$ of \cref{edge bound on theta-free graphs}, we make three pertinent remarks.
Firstly, although $2n-2$ is upper bounded by $3n-6$ (for $n \geqslant 4$), \tf{} graphs are not necessarily planar; consider the graph $\tftwelve$ shown in \cref{fig:barriertf}.
Secondly, it follows immediately that every \tf{} graph has a vertex of degree at most three.

\begin{corollary}
	Every matching covered graph, with minimum degree $\delta \geqslant 4$, is \tb{}. 
	\qed
\end{corollary}

Finally, there does not exist a constant upper bound on the maximum degree of \tf{} graphs, as stated below and proved in \cref{subsec:sparsity and bicriticality}.

\begin{theorem}\label{for any k there is a tf graph with a vertex of degree k}
	For any positive integer $k$, there exists a \tf{} graph with a vertex of degree $k$.
\end{theorem}

Recall that a \mcg{} has zero bricks if and only if it is bipartite; on the other hand, it has zero braces if and only if it is bicritical and has order four or more.
In this sense, bicritical graphs (of order four or more) may be viewed as those \mcg{}s that are `farthest' from being bipartite.
It turns out that the members of $\mathcal{T}$ are essentially all of the \tf{} bicritical graphs as stated below.

\begin{restatable}{theorem}{tfBicritical}
	{\sc[\tf{} Bicritical Graphs]}\label{characterization of theta-free-bicritical graphs} \newline
	A \tf{} graph $G$ is bicritical if and only if $G \in \{K_2\} \cup \mathcal{T}$.
\end{restatable}

A proof appears in \cref{subsec:sparsity and bicriticality}.
In contrast to the above, we now proceed to discuss a class of nonbipartite \mcg{}s that are `closest' to being bipartite.

			\subsubsection{\tf{} near-bipartite graphs}
				\label{subsec:tf near-bipartite intro}
				We begin by recalling \cref{the only bipartite theta-free mcgs are k2 and even cycles}.
In the same spirit, the class of \tf{} `near-bipartite graphs' turns out to be very restrictive.
In order to define them, and to illustrate their significance, we briefly discuss the ear decomposition theory from \cite[Chapter 11]{lumu24}.

By the \emph{length} of a path, we mean the number of its edges.
A \emph{single ear} of a graph $G$ is an odd path $P$ each of whose internal vertices (if any) has degree two in $G$;
the graph $G-P$ denotes the one obtained from $G$ by deleting all edges and internal vertices of $P$, and we say that it is obtained from $G$ by \emph{deleting the single ear $P$}.
A single ear $P$ of a \mcg{} $G$ is \emph{removable} if $G-P$ is also matching covered.
Hetyei proved that every bipartite \mcg{}, except $K_2$, has a removable single ear; see \cite[Theorem 3.14]{lumu24}.
However, this does not hold for nonbipartite \mcg{}s; consider $K_4$ or $\overline{C_6}$.

A \emph{double ear} of a \mcg{} $G$ is a pair of vertex-disjoint single ears $\{P,Q\}$;
furthermore, it is \emph{removable} if neither $P$ nor $Q$ is removable, but $G-P-Q$ is \mc{}, and we say that $G-P-Q$ is obtained from $G$ by \emph{deleting the double ear $\{P,Q\}$}.
For instance, observe that each perfect matching of $K_4$ is a removable double ear.
A \emph{removable ear} is either a single or a double ear that is removable.
Now, we state the Ear Decomposition Theorem \cite[Theorem 11.8]{lumu24} due to Lov\'asz and Plummer.

\begin{theorem}{\sc[Existence of Removable Ears]}
	\label{existence of removable ears}
	\newline
	Every \mcg, except $K_2$, has a removable ear.
\end{theorem}

Thus, for any \mcg{} $G$, there exists a sequence of \mc{} subgraphs $(K_2=G_0,G_1,\dots,G_{r-1},G_r=G)$ such that, for each $1 \leqslant i \leqslant r$, the graph $G_{i-1}$ is obtained from $G_i$ by deleting a removable ear of $G_i$, and such a sequence is called an \emph{ear decomposition} of $G$; we invite the reader to observe that $G_{i-1}$ is a conformal subgraph of $G_i$, and thus of every subsequent element.
It follows from the aforementioned result of Hetyei that if $G$ is nonbipartite, then the first nonbipartite subgraph, in any ear decomposition, has a removable double ear whose deletion results in a bipartite graph.
This brings us to the definiton of \emph{near-bipartite} graphs: any (nonbipartite) \mcg{} that has a removable double ear whose deletion results in a bipartite graph.
For example, $K_4$ and $\overline{C_6}$ are near-bipartite, whereas $\petersen$ is not.

The above discussion sheds light on the relation between ear decompositions and conformal \mc{} subgraphs.
In particular, the second element in any ear decomposition is a conformal cycle.
A \mcg{} $G$ is \emph{cycle-extendable} if each of its even cycles is conformal, or equivalently, using a stronger version of \cref{existence of removable ears} (\cite[Theorem 11.7]{lumu24}), each of its even cycles is the second element in some ear decomposition of $G$.
Recently, Dalwadi, Pause, Diwan and Kothari \cite{dpdk25} obtained a characterization of planar cycle-extendable graphs, and we refer the reader to their work for further context.
We are now ready to state our characterizations of near-bipartite \tf{} graphs; a proof appears in \cref{subsec:tf near-bipartite}.

\begin{restatable}{theorem}{tfNearBipartite}{\sc[Near-bipartite \tf{} Graphs]}\label{tf near-bipartite} \newline
	For a \tf{} graph $G$, the following are equivalent:
	\begin{enumerate}[label=(\roman*)]
		\item $G$ is either bipartite or near-bipartite,
		\item either $b(G) = 0$, or $b(G)=1$ and $K_4$ is the unique brick of $G$,
		\item $G$ is cycle-extendable, and
		\item $G$ is either $K_2$, an (even) cycle, or a bisubdivision of $K_4$.
	\end{enumerate}
\end{restatable}

We conclude this section by stating a stronger version of Lov\'asz's Theorem (\ref{k4-c6bar theorem}): every nonbipartite \mcg{} admits an ear decomposition wherein either the second element is a bisubdivision of $K_4$, or the third one is a bisubdivision of $\overline{C_6}$.

			\subsubsection{Pfaffian graphs}
				\label{subsubsec:tf pfaffian intro}
				We remark that it suffices to restrict attention to \mcg{}s, even though some of the following may be generalized to all graphs.
Famously, Valiant \cite{vali79} showed that counting the number of perfect matchings in a given graph is \mbox{$\#\mathsf{P}$-complete} even for bipartite graphs.
However, certain classes, including Pfaffian graphs defined below, admit a polynomial time algorithm for the same.
The interested reader may refer to \cite{thom06} and \cite[Part III]{lumu24} for additional context.

With respect to an orientation $D$ of a graph $G$, an even cycle of $G$ is \emph{oddly oriented} if the number of its edges that agree with either prescribed sense of traversal is odd.
For instance, in the graph shown in \cref{fig:Pfaffian orientation example}, the $4$-cycle $C$ (highlighted in blue) is oddly oriented, whereas the $6$-cycle~$D$ (shown in thick pink lines) is not.
% For instance, in the graph shown in \cref{fig:Pfaffian orientation example}, the $4$-cycle $C$ is oddly oriented, whereas the $6$-cycle $D$ (shown in pink) is not.
An orientation of a \mcg{} is \emph{Pfaffian} if each conformal cycle is oddly oriented;
\cref{fig:Pfaffian orientation example} shows an example.

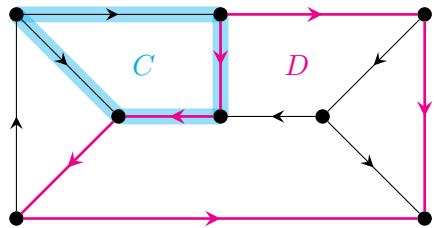
\fig{A Pfaffian orientation of the bicorn graph}{fig:Pfaffian orientation example}{\begin{tikzpicture}[scale=1.35]

\tikzstyle{arrow}=[postaction={decorate,decoration={markings,mark=at position 0.5 with {\arrow[scale=1.5,>=stealth]{>}}}}]
\tikzstyle{arrowcolor}=[line width=1pt,postaction={decorate,decoration={markings,mark=at position 0.5 with {\arrow[scale=1.25,>=stealth]{>}}}}]

\draw[line width=6pt, opacity=0.4, cyan] (-2,1)--(-1,0);
\draw[line width=6pt, opacity=0.4, cyan] (-2,1)--(0,1);
\draw[line width=6pt, opacity=0.4, cyan] (0,0)--(-1,0);
\draw[line width=6pt, opacity=0.4, cyan] (0,1)--(0,0);

\node[vtx] (l1) at (-2,1) {};
\node[vtx] (l2) at (-1,0) {};
\node[vtx] (l3) at (-2,-1) {};
\node[vtx] (r1) at (2,1) {};
\node[vtx] (r2) at (1,0) {};
\node[vtx] (r3) at (2,-1) {};
\node[vtx] (m2) at (0,0) {};
\node[vtx] (m1) at (0,1) {};

\draw[arrow] (l3)--(l1);
\draw[arrow] (l1)--(l2);
\draw[arrowcolor,magenta] (l2)--(l3);

\draw[arrow] (l1)--(m1);
\draw[arrowcolor,magenta] (m1)--(r1);
\draw[arrow] (r2)--(m2);
\draw[arrowcolor,magenta] (m2)--(l2);
\draw[arrowcolor,magenta] (l3)--(r3);

\draw[arrowcolor,magenta] (m1) -- (m2);

\draw[arrowcolor,magenta] (r1)--(r3);
\draw[arrow] (r1)--(r2);
\draw[arrow] (r2)--(r3);

\node[magenta] at (0.75,0.5) {$D$};
\node[cyan] at (-0.75,0.5) {$C$};

\end{tikzpicture}}

If $D$ is a Pfaffian orientation of a \mcg{} $G$, (the absolute value of) the determinant of the adjacency matrix of~$D$ is the square of the number of perfect matchings of $G$, and thus the latter may be computed in polynomial time (as alluded to earlier).
This inspires the following definition and decision problems;
a \mcg{} is \emph{Pfaffian} if it admits a Pfaffian orientation.

\begin{decision-problem}
	\label{dp:is D Pfaffian}
	Given an orientation of a \mcg{}, decide whether it is Pfaffian.
\end{decision-problem}

\begin{decision-problem}
	\label{dp:is G Pfaffian}
	Given a matching covered graph, decide whether it is Pfaffian.
\end{decision-problem}

Carvalho, Lucchesi and Murty \cite[Theorem 3.9]{clm05a} used ear decompositions to devise a polynomial time algorithm to compute a Pfaffian orientation, provided the input graph is Pfaffian.
Using this, they deduced that \cref{dp:is G Pfaffian,dp:is D Pfaffian} are computationally equivalent.
They also showed that \cref{dp:is G Pfaffian} is in \coNP{}.
These results were first established by Vazirani and Yannakakis \cite{vaya89} using linear algebraic techniques.
Neither \coNP{} certificate is in terms of conformal minors; we refer the reader to \cite[Chapter 20]{lumu24}.

Kasteleyn \cite{kast63} showed that every planar graph is Pfaffian.
The converse does not hold; for instance, the Heawood graph, shown in \cref{fig:heawood}, is Pfaffian.
It is well known that $K_{3,3}$ is the smallest non-Pfaffian graph (see \cite[Lemma 19.10]{lumu24}).
Little \cite{litt75} established that $K_{3,3}$ is the only minimally non-Pfaffian graph amongst the bipartite ones (see \cref{characterization of bipartite Pfaffian graphs})  --- where \emph{minimally non-Pfaffian} means that the graph is non-Pfaffian, whereas each of its proper conformal minors is Pfaffian.

The Petersen graph $\petersen$ is non-Pfaffian~\cite[Lemma 19.11]{lumu24}, and it is the smallest one amongst \tf{} graphs; the graph $\tftwelve$, shown in \cref{fig:barriertf}, is the next smallest one except for the unique bisubdivision of $\petersen$ of order twelve.
The equivalence of $(i)$ and $(iii)$, in our characterization of Pfaffian $\tf{}$ graphs stated below, shows that these are the only minimally non-Pfaffian \tf{} graphs.

\begin{restatable}{theorem}{tfPfaffian}{\sc[Characterization of Pfaffian \tf{} Graphs]}\label{theta-free Pfaffian characterization} \newline
	For a \tf{} graph $G$, the following are equivalent:
	\begin{enumerate}[label=(\roman*)]
		\item $G$ is Pfaffian,
		\item $G$ does not contain either of $K_{3,3}$ and the Petersen graph $\petersen$ as an \sminor,
		\item $G$ does not contain either of $\tftwelve$ and $\petersen$ as a conformal minor, that is, $G$ is \mbox{$\tftwelve$-free} and \mbox{$\petersen$-free}, and
		\item each brace of $G$ is Pfaffian, and $\petersen$ is not a brick of $G$.
	\end{enumerate}
\end{restatable}

The `$S$-minor' containment notion is defined in \cref{subsec:tf Pfaffian}.
In the case of \tf{} graphs, the equivalence of $(i)$ and $(iii)$ yields a \coNP{} certificate for \cref{dp:is G Pfaffian} in terms of conformal minors, \`a la Little's result \cite{litt75}.
The above theorem, which is proved in \cref{subsec:tf Pfaffian}, is reminiscent of the following.
% The equivalence of $(i)$ and $(iv)$, in the above, implies that \cref{dp:is G Pfaffian} is in \PP{} for \tf{} graphs.
% in order to prove it (in \cref{subsec:tf Pfaffian}), we make use of the aforementioned result of Little (1975), stated below.

\begin{theorem}{\sc[Characterization of Pfaffian Bipartite Graphs]}\newline
	\label{characterization of bipartite Pfaffian graphs}
	For a bipartite matching covered graph $H$, the following are equivalent:
	\begin{enumerate}[label=(\roman*)]
		\item $H$ is Pfaffian,
		\item $H$ does not contain $K_{3,3}$ as an \sminor,
		\item $H$ does not contain $K_{3,3}$ as a conformal minor, that is, $H$ is \mbox{$K_{3,3}$-free}, and
		\item each brace of $H$ is Pfaffian.
	\end{enumerate}
\end{theorem}

The equivalence of $(i)$ and $(iv)$ is a special case of a result first proved by Little and Rendl~\cite{lire91}.
On the other hand, the equivalence of $(i)$ and $(iii)$ is due to the aforementioned pioneering work of Little.
However, it was not until more than two decades later that Robertson, Seymour and Thomas \cite{rst99}, and independently McCuaig \cite{mccu04}, established an \NP{}-characterization of Pfaffian bipartite graphs; see \cite{mrst97}.
Their breakthrough result also yields a polynomial time algorithm, thus placing \cref{dp:is G Pfaffian} in \PP{} for bipartite graphs.
This, along with the equivalence of~$(i)$~and~$(iv)$ stated in \cref{theta-free Pfaffian characterization}, and the fact that bricks and braces of a \mcg{} may be computed efficiently, yields a polynomial time algorithm for recognizing Pfaffian \tf{} graphs; ergo, placing \cref{dp:is G Pfaffian} in \PP{} for \tf{} graphs.

In the bipartite case, one may impose a stronger condition to define a proper subset of Pfaffian graphs.
An equivalent formulation of this condition led us to a couple of open problems, and to a serendipitous connection between one of these and \tf{}ness; see \cref{a cubic graph is theta free iff each conformal cycle is of length 0 mod 4}.
We discuss these next.

			\subsubsection{Lengths of conformal cycles modulo 4}
				The \emph{natural orientation} $D$ of a bipartite \mcg{}~$H[A,B]$ is the one obtained by directing each edge from its end in~$A$ to its end in~$B$.
One may now ask whether this orientation is Pfaffian or not.
Due to the equivalence of $(i)$ and $(ii)$ stated below, $H$ is said to be \emph{det-extremal} if $D$ is Pfaffian; see \cite{fjls03}.

\begin{proposition}
	For a bipartite \mcg{}, the following are equivalent:
	\begin{enumerate}[itemsep=0pt,label=(\roman*)]
		\item its natural orientation is Pfaffian,
		\item the determinant and the permanent, of its bipartite adjacency matrix, have equal absolute values,
		\item each of its conformal cycles is of length $2 \pmod{4}$.
	\end{enumerate} 
\end{proposition}

The reader may observe the equivalence of $(i)$ and $(iii)$.
Inspired by the above, we propose the following open problem.

\begin{problem}
	\label{problem:cc 2 mod 4}
	Characterize \mcg{}s each of whose conformal cycles is of length $2 \pmod{4}$.
\end{problem}

Each conformal cycle of the Heawood graph, shown in \cref{fig:heawood}, is of length $2 \pmod{4}$;
% McCuaig~\cite{mccu98} proved that amongst cubic braces of order four or more, it is the only such graph, and thus deduced the following.
McCuaig~\cite{mccu98} proved that it is the only such brace (of order four or more), and also deduced the following.

\begin{theorem}
	For a \con3 bipartite cubic graph, each of its conformal cycles is of length $2 \pmod{4}$ if and only if each of its braces is (isomorphic to) the Heawood graph.
\end{theorem}

\fig{The Heawood graph}{fig:heawood}{\begin{tikzpicture}[scale=1]

\tikzmath{
	\s = 360/7+90;
	\a = 1;
	\b = 2.1;
}

% \foreach \i in {1,2,...,7}{
% 	\node[vtx-white] (u\i) at (\s+\i*360/7:\a) {};
% 	\node[vtx] (v\i) at (\s+\i*360/7:\b) {};
% 	\draw (u\i) -- (v\i);
% }

% \foreach \i in {1,2,...,7}{
% 	\pgfmathtruncatemacro{\j}{\i+2}
% 	\pgfmathtruncatemacro{\k}{\i-1}
% 	\pgfmathtruncatemacro{\l}{\i+1}

% 	\ifnum\j>7
% 		\pgfmathtruncatemacro{\j}{\j-7}
% 	\fi
% 	\ifnum\l>7
% 		\pgfmathtruncatemacro{\l}{\l-7}
% 	\fi
% 	\ifnum\k<1
% 		\pgfmathtruncatemacro{\k}{7}
% 	\fi

% 	% \draw (u\i) -- (v\j);
% 	% \draw (u\i) -- (v\k);
	
% 	\draw (u\i) to[bend left=20] (v\k);
% 	\draw (u\i) to[bend right=30] (v\j);
% }

\begin{scope}[shift={(-6.5,0)}]

	\tikzmath{
		\s = 360/14;
		\a = 1;
		\b = 2.1;
	}

	\foreach \i in {1,...,14}{
		\pgfmathtruncatemacro{\j}{\i-1}
		\ifodd\i
			\node[vtx-white] (\i) at (\s+\i*360/14:\b) {};
			% \node at (\s+\i*360/14:2.25) {$\i$};
		\else
			\node[vtx] (\i) at (\s+\i*360/14:\b) {};
			% \node at (\s+\i*360/14:2.25) {$\i$};
		\fi
		\ifnum\j>0;
			\draw (\j) -- (\i);
		\fi
	}
	\draw (14) -- (1);

	\draw (1) -- (10);
	\draw (2) -- (7);
	\draw (3) -- (12);
	\draw (4) -- (9);
	\draw (5) -- (14);
	\draw (6) -- (11);
	\draw (8) -- (13);

\end{scope}

\end{tikzpicture}}

To the best of our knowledge, the \con2 bipartite cubic case of \cref{problem:cc 2 mod 4} is still open.
It is only natural that we also put forth the analogous open problem stated below.

\begin{problem}
	\label{problem:cc 0 mod 4}
	Characterize \mcg{}s each of whose conformal cycles is of length $0 \pmod{4}$.
\end{problem}

It turns out that \cref{problem:cc 0 mod 4,problem:cc 2 mod 4} are somewhat related to \cref{theta-k4 theorem}, as we explain next.
Let~$H$ be any bisubdivision of $\theta$.
Note that $\theta$ has precisely three even cycles, each of which is the complement of a perfect matching.
These correspond to (conformal) cycles of $H$; let their lengths be $\ell_1,\ell_2$ and $\ell_3$.
Note that each edge of $H$ is counted precisely twice in $\ell:=\ell_1+\ell_2+\ell_3$. Since $H$ is a bisubdivision of~$\theta$, we infer that $\ell \equiv 2 \pmod{4}$.
Consequently, at least one of $\ell_1,\ell_2$ and $\ell_3$ is congruent to $2 \pmod{4}$.
This, combined with transitivity of conformality, proves the first statement of the following;
whereas the second one follows by analogous arguments for $K_4$, in which case $\ell \equiv 0 \pmod{4}$.

\begin{proposition}\label{theta and K4 vs parity of conformal cycles}
	Every \tb{} \mcg{} has a conformal cycle of length $2 \pmod{4}$, whereas
	every $K_4$-based \mcg{} has a conformal cycle of length $0 \pmod{4}$. \qed
\end{proposition}

Let us now see an easy, yet pertinent, application of the above.

\begin{proposition}\label{Petersen is theta free}
	The Petersen graph $\mathbb{P}$ is \tf{}.
\end{proposition}
\begin{proof}
	It is prudent to exploit the symmetries of the Petersen graph.
	Since $\mathbb{P}$ has girth five and is non-hamiltonian, its only even cycles are $6$-cycles and $8$-cycles.
	If $C$ is any $8$-cycle, then $\mathbb{P}-V(C)$ is $K_2$; whereas if $C$ is any $6$-cycle, then $\mathbb{P}-V(C)$ is $K_{1,3}$.
	Thus, the conformal cycles of $\mathbb{P}$ are precisely its $8$-cycles.
	Consequently, by \cref{theta and K4 vs parity of conformal cycles}, $\mathbb{P}$ is \tf.
\end{proof}

For either clause of \cref{theta and K4 vs parity of conformal cycles}, the converse does not hold.
For instance, $\overline{C_6}$ is a (\con3 cubic) $K_4$-free graph that has conformal cycles of length four;
whereas the \tf{} graph $T_6$, shown in \cref{fig:T6ill}, has conformal cycles of length six.
Interestingly, in the case of \cc2s, the converse of the first clause of \cref{theta and K4 vs parity of conformal cycles} holds; a proof appears in \cref{subsec:tf cubic}.

\begin{restatable}{theorem}{tfCubic}
	\label{a cubic graph is theta free iff each conformal cycle is of length 0 mod 4}
	{\sc[Cubic \tf{} Graphs vs Lengths of Conformal Cycles modulo 4]} \newline
	For any $2$-connected cubic graph $G$, each of its conformal cycles is of length $0 \pmod{4}$ if and only if $G$ is \tf{}.
\end{restatable}

The above is the serendipitous discovery that we alluded to earlier, and was in fact first envisaged by the third author.
\cref{a cubic graph is theta free iff each conformal cycle is of length 0 mod 4} solves \cref{problem:cc 0 mod 4} for the case of \cc2s, and yields a polynomial time algorithm for the corresponding decision problem.
This is particularly fascinating given that \cref{problem:cc 2 mod 4} is open even in the case of \con2 bipartite cubic graphs.
Speaking of cubic graphs, one is also tempted to discuss edge-colorings, and so we do.

			\subsubsection{3-edge-colorability of cubic graphs}
				\label{subsubsec:tight cuts}
				It is well known that the edge chromatic number of any cubic graph is either three or four.
This may be viewed as a special case of the Vizing-Gupta Theorem \cite[Theorem 17.4]{bomu08}.
Note that a connected cubic graph, that is not \con2, requires four colors;
this leads us to the following decision problem.

\begin{decision-problem}
	\label{dp:cubic-3ec}
	Given a \cc2, decide whether it is $3$-edge-colorable.
\end{decision-problem}

The above is clearly in \NP{}, but is not known to be in \coNP{};
Holyer \cite{holy81}, and independently Galil and Leven \cite{gale83}, showed that it is in fact \mbox{$\mathsf{NP}$-hard}.
The famous result of Tait \cite[Theorem 11.4]{bomu08} established the following equivalent version of the Four Color Theorem: every planar \cc2 is \mbox{$3$-edge-colorable}.
In light of this, recall from \cref{subsubsec:sparsity and bicriticality intro} that \tf{} graphs, including cubic ones, need not be planar;
one such example is the Petersen graph, and it is famously known to be the smallest \cc2 that is not $3$-edge-colorable.
In \cref{subsec:tf cubic 3ec}, using our main result, we deduce the following characterization of \tf{} cubic graphs that are \mbox{$3$-edge-colorable}.

\begin{restatable}{theorem}{tfCubicThreeEdgeColorability}
	\label{3-e-c cubic theta free graphs}
	{\sc[$3$-edge-colorability of \tf{} Cubic Graphs]}\newline
	For a \tf{} cubic graph $G$, the following are equivalent:
	\begin{enumerate}[label=(\roman*)]
		\item $G$ is $3$-edge-colorable,
		\item each brick of $G$ is isomorphic to $K_4$,
		\item the Petersen graph is not a brick of $G$,
		\item $G$ does not contain the Petersen graph as a conformal minor, that is, $G$ is \mbox{$\petersen$-free}, and
		\item $G$ does not contain the Petersen graph as an \sminor{}.
	\end{enumerate}
\end{restatable}

In the case of \tf{} graphs, the above places \cref{dp:cubic-3ec} in \coNP{}; it also yields a polynomial time algorithm, and thus places it in \PP{} as well.

		\subsection{Organization of this paper}
			In \cref{sec:characterizing theta-free graphs}, we prove the Main Theorem (\ref{theta-free characterization - inductive version}).
In \cref{sec:consequences of main theorem}, we elaborate on, and deduce all of the consequences in the same order as they appear in \cref{subsec:consequences}.

	\section{Proof of the Main Theorem}
		\label{sec:characterizing theta-free graphs}
		We begin this section by stating \cref{theta-free characterization - inductive version} from a generation view point.
To this end, we define an operation which may be viewed as the inverse of shriniking a shore of a cut.

Let $v_1$ and $v_2$ be vertices of disjoint graphs $G_1$ and $G_2$, respectively, and let $\pi:\partial_{G_1}(v_1) \rightarrow \partial_{G_2}(v_2)$ denote a bijection.
Let $G$ be the graph obtained by adding to $G_1-v_1 \,\cupdot\, G_2-v_2$ the following new edges:
for each edge $e$ in $\partial_{G_1}(v_1)$, add a corresponding edge joining the end of $e$ distinct from~$v_1$ and the end of $\pi(e)$ distinct from~$v_2$.
We say that $G$ is obtained by \emph{splicing $G_1$ and $G_2$ at the vertices $v_1$ and $v_2$ with respect to $\pi$}, and we denote this as $G:=(G_1 \odot G_2)_{v_1,v_2,\pi}$; if the choices of $v_1,v_2$ and $\pi$ are either irrelevant or clear from context, we simplify the notation to $G:=G_1 \odot G_2$.
The cut $C:=\partial_G(V(G_1)-v_1)$ is called the corresponding \emph{splicing cut}; observe that $G_1$ and $G_2$ are the two $C$-contractions of $G$.
It is easy to see that cut contractions are unique but splicing two graphs may result in distinct graphs, based on the choice of $v_1,v_2$ and~$\pi$.
The reader may also verify that any splicing of two \mcg{}s yields a \mcg{}.

We are now ready to state the generation view point of \cref{theta-free characterization - inductive version}.

\begin{theorem}{\sc[Characterization of \tf{} Graphs - Generation Version]}\label{theta-free characterization - generation version}
	\newline
	A matching covered graph~$G$ is \tf{} if and only if at least one of the following holds: 
	\begin{enumerate}[label=(\roman*)]
		\item either $G$ is $K_2, C_2, K_4$ or the Petersen graph $\mathbb{P}$, or
		\item $G$ is obtained from a bipartite \mcg{} $H[A,B]$ by splicing at each vertex $a_i \in A$ with some \tf{} graph $G_i$,~or
		\item $G$ is a $K_2$-sum of two \tf{} graphs.
	\end{enumerate}
\end{theorem}

Next, we discuss a couple of observations that we shall find useful in proving \cref{J-freeness across a barrier where J is bicritical cubic,J-freeness across a 2-separation where J is 3-connected cubic}.
In particular, we consider the intersection of induced paths of a conformal subgraph with a tight cut, and state some of the terminology that was introduced by Kothari and Murty \cite{km16}.

		\subsection{$C$-crossing paths}
			\label{subsec:C-crossing paths}
			For a subdivision $H$ of a graph~$J$ with $\delta(J)\geqslant 3$, we refer to the vertices of $H$ of degree three or more as its {\it branch vertices}, and to its remaining vertices as {\it subdivision vertices}.
Note that each branch vertex corresponds to a unique vertex of $J$ and vice versa; we use the same label to refer to both.
One may associate with each edge $uv\in E(J)$ a unique path $P_{uv}$ in $H$, each of whose internal vertices (possibly none) is a subdivision vertex of $H$.
We remark that this association is uniquely defined whenever $J$ is simple; however, since we are particularly interested in $J=\theta$, we do not assume its uniqueness, but rather fix any particular association.

In case $H$ is a bisubdivision of $J$, we may assign a parity to each edge of $P_{uv}$  (and thus of $H$) as described below.
Let $P_{uv}:= w_1w_2\dots w_{2k}$ where $w_1:=u$ and $w_{2k}:=v$.
We say that an edge $w_iw_{i+1}$ is of {\it odd parity} if $i$ is odd; otherwise, it is of {\it even parity}.
Note that the parity of an edge is independent of the order in which the path is traversed.
In particular, the first and last edges (possibly not distinct) of any such path, say $P_{uv}$, are both of odd parity.

Now, we assume $J$ to be matching covered and relate perfect matchings of $J$ with those of~$H$.
Let $M_J$ denote any perfect matching of $J$.
We may obtain a perfect matching $M_H$ of $H$ as follows.
For each edge $uv \in E(J)$, if $uv \in M_J$ then include the odd edges of $P_{uv}$ into $M_H$; otherwise, include the even edges of $P_{uv}$.
It is easy to see that $M_H$ is indeed a perfect matching of $H$.
We invite the reader to observe that every perfect matching of $H$ arises in this manner, and that $M_J \mapsto M_H$ is a bijection from the collection of perfect matchings of $J$ to that of $H$.
This also proves \cref{every bisubdivision of a mcg is mc} discussed earlier.

Now, suppose that $H$ is a conformal subgraph of a \mcg{} $G$, and let $C$ denote any cut of $G$.
For an edge $uv$ of $J$, we say that the path $P_{uv}$ (of $H$) is a {\it $C$-crossing path} if at least one edge of $P_{uv}$ lies in $C$, and we call such an edge a {\it $C$-crossing edge}.
If there is one such edge, we say that $P_{uv}$ crosses $C$ once; if there are two such edges, we say that $P_{uv}$ crosses $C$ twice; and so on.
We now state some observations due to Kothari and Murty \cite{km16}.
Although they imposed stronger conditions on $J$, their proof goes through exactly as is with the weaker assumptions stated below.

\begin{proposition}\label{observations on C-crossing paths}
	Let $G$ be a \mcg{} that contains a conformal bisubdivision $H$ of~$J$, where $J$ is a \mcg{} with \mbox{$\delta(J)\geqslant 3$}.
	Then, for any tight cut~$C$ of $G$, each of the following holds.
	\begin{enumerate}[label=(\roman*)]
		\item For a $C$-crossing path $P_{uv}$, any two $C$-crossing edges of $P_{uv}$ must be of opposite parity.
		Thus, $P_{uv}$ crosses $C$ at most twice.
		\item If a $C$-crossing path $P_{uv}$ crosses $C$ twice, then there are no other $C$-crossing paths.
		\item If $P_{uv}$ and $P_{uw}$ are two $C$-crossing paths, then each of them must cross $C$ in an edge of odd parity.
	\end{enumerate}
\end{proposition}

When convenient, we shall view subgraphs as edge sets and vice versa.
We now prove an easy consequence of the above that we shall find useful in proving \cref{J-freeness across a barrier where J is bicritical cubic}.

\begin{corollary} \label{if a tight cut shore contains at most one branch vertex then J-basedness is preserved}
	Let $G$ be a \mcg{} that contains a conformal bisubdivision $H$ of~$J$, where $J$ is a \mcg{} with \mbox{$\delta(J)\geqslant 3$}, and let $C:=\partial(X)$ be a tight cut of $G$.
	If the shore $\overline{X}$ contains at most one branch vertex of $H$, then $G/\overline{X}$ is $J$-based.
\end{corollary}
\begin{proof}
	Let $G':=G/\overline{X}$, let $M_H$ denote a perfect matching of $H$, and let $M$ be an extension of $M_H$ to a perfect matching of $G$.
	We consider a couple of cases and make some observations pertaining to the cardinalities of $M_H \cap C$ and $(M-M_H) \cap C$.

	First, suppose that all branch vertices lie in $X$.
	Then by \cref{observations on C-crossing paths} $(i)$ and $(ii)$, $|\partial(X)\cap E(H)|\in\{0,2\}$.
	Note that in the first case $|M_H\cap C|=0$ and $|(M-M_H)\cap C|=1$; whereas, in the second case, both edges belong to the same $C$-crossing path and are of opposite parities, and thus $|M_H\cap C|=1$ and $|(M-M_H)\cap C|=0$.
	Now, suppose that there is a branch vertex $z$ in $\overline{X}$.
	Since $\delta(J)\geqslant 3$, by \cref{observations on C-crossing paths} $(iii)$, each $C$-crossing path emanating from $z$ meets $C$ in an edge of odd parity.
	Consequently, $|M_H\cap C|=1$ and $|(M-M_H)\cap C|=0$.

	Observe that in all the cases discussed above, the subgraph $H':=G'\cap E(H)$ is a bisubdivision of $J$, and that $(M-M_H)\cap E(G')$ is a perfect matching of $G'-V(H')$.
	This completes the proof of \cref{if a tight cut shore contains at most one branch vertex then J-basedness is preserved}.
\end{proof}

In the next three sections, we prove \cref{J-freeness across a barrier where J is bicritical cubic,J-freeness across a 2-separation where J is 3-connected cubic,K4 and Petersen are the only theta-free bricks}, respectively.

		\subsection{Dealing with barriers: a proof of \cref{J-freeness across a barrier where J is bicritical cubic}}
			\label{sec:barrier}
			We start by proving a couple of technical results.

\begin{lemma}\label{unique branch vertices of H in distinct components of G-B are nonadjacent in J}
	Let $G$ be a \mcg{} that contains a conformal bisubdivision $H$ of~$J$, where $J$ is a \mcg{} with \mbox{$\delta(J)\geqslant 3$}, and let $B$ denote a barrier of $G$.
	If $G-B$ has distinct components, say $K$~and~$L$, each of which contains exactly one branch vertex of~$H$, say~$u$~and~$v$, respectively, then $u$~and~$v$ are nonadjacent in $J$.  
\end{lemma}
\begin{proof}
	Suppose to the contrary that $u$~and~$v$ are adjacent in $J$, and let $P_{uv}$ denote the path of $H$ corresponding to an edge joining $u$~and~$v$ in $J$.
	Since $\delta(J)\geqslant 3$, there are three or more \mbox{$\partial(K)$-crossing} paths emanating from $u$; by \cref{observations on C-crossing paths} $(iii)$, each of them crosses $\partial(K)$ in an edge of odd parity. In particular, $P_{uv}$ meets $\partial(K)$ in an edge of odd parity, say $e_1$. By an analogous argument,~$P_{uv}$ meets $\partial(L)$ in an edge of odd parity, say $e_2$. Consequently, the subpath $e_1P_{uv}e_2$ is an odd path. We intend to arrive at a contradiction to this.

	Note that if $P_{uv}$ meets a component $F$ of $G-B$ that is distinct from $K$ and $L$, then by \cref{observations on C-crossing paths} $(i)$, $P_{uv}$ has precisely two $\partial(F)$-crossing edges and they are of opposite parity; thus, the subpath $P_{uv}\cap F$ is even. Using this observation, and the fact that $B$ is a stable set, we conclude that the subpath $e_1P_{uv}e_2$ is an even path; this contradicts what we have already established above.  
\end{proof}

We now use the above result to prove a generalization of \cref{any claw of a bip mcg is part of a conformal bisubdivision of theta}.

\begin{proposition}
	\label{any claw centered at an isolated vertex of G-B participates in a conformal bisubdivision of theta}
	Let $B$ denote a barrier of a matching covered graph $G$, and let $z$ denote an isolated vertex of $G-B$.
	Then, any three edges in $\partial(z)$ participate in a conformal bisubdivision of~$\theta$.
\end{proposition}
\begin{proof}
	Let $K$ denote a set of three edges incident at~$z$.
	We proceed by induction on the order of~$G$.
	If each (odd) component of $G-B$ is trivial, then $G$ is bipartite, the desired conclusion holds by \cref{any claw of a bip mcg is part of a conformal bisubdivision of theta}.
	Now let $L$ be a nontrivial odd component of $G-B$.
	Let $U:=V(L)$, let $L_1:=G/U\rightarrow u$, and let $L_2:=G/\overline{U}\rightarrow \overline{u}$.
	Observe that $B$ is a barrier of $L_1$, and that $z$ is an isolated vertex of $L_1-B$.
	By the induction hypothesis, $K$ participates in a conformal bisubdivision of~$\theta$, say $H_1$, in $L_1$; let $M_1$ be a perfect matching of $L_1-V(H_1)$.
	By \cref{unique branch vertices of H in distinct components of G-B are nonadjacent in J}, the branch vertex of $H_1$, distinct from $z$, lies in $B$; in particular, $u$~is not a branch vertex of $H_1$.
	
	Now, let $F:=E(H_1)\cap\partial(U)$.
	It follows from the above and \cref{observations on C-crossing paths} $(i)$ and $(ii)$ that $|F|\in\{0,2\}$.
	If $|F|=0$, let $H_2$ denote the null subgraph of $L_2$ and let $M_2$ denote a perfect matching of $L_2$ that contains the unique edge in $M_1\cap\partial(U)$.
	On the other hand, if $|F|=2$, we invoke \cref{any two adjacent edges participate in a conformal cycle and C2 is the only nonsimple theta-free mcg}, and let $H_2$ denote a conformal cycle in $L_2$ containing $F$, and let $M_2$ denote a perfect matching of $L_2-V(H_2)$.
	Observe that, in each case, $H:=H_1\cup H_2$ is a bisubdivision of~$\theta$ in $G$, and that $M:=M_1\cup M_2$ is a perfect matching of $G-V(H)$.
\end{proof}

We shall now prove the forward implication of \cref{J-freeness across a barrier where J is bicritical cubic} with the weaker hypothesis that~$J$ is a \mcg{} with maximum degree $\Delta(J) \leqslant 3$.

\begin{theorem}\label{G is J-based if there is a B-fragment that is J-based where J has max degree 3}
	Let $B$ denote a barrier of a \mcg{} $G$, and let $J$ be a \mcg{} with $\Delta(J)\leqslant 3$. If there exists \mbox{$B$-fragment} that is $J$-based, then $G$ is also $J$-based.
\end{theorem}
\begin{proof}
	Since every \mcg{} is $K_2$-based, we may assume that $J\neq K_2$; henceforth, $\delta(J)\geqslant 2$.
	Let $H_1$ denote a conformal bisubdivision of $J$ in a $B$-fragment $L_1$ corresponding to an odd component $L$ of $G-B$, and let $M_1$ denote a perfect matching of $L_1-V(H_1)$.
	Let $F:=E(H_1)\cap\partial_{L_1}(L)$; since \mbox{$2\leqslant\delta(J)\leqslant\Delta(J)\leqslant 3$}, we infer that $|F|\in\{0,2,3\}$.
	Let $L_2$ denote the other \mbox{$\partial(L)$-contraction} of $G$; that is $L_2:=G/V(L)$.
	Depending on the cardinality of $F$, we will define a conformal subgraph $H_2$ in $L_2$, and a perfect matching $M_2$ of $L_2-V(H_2)$.

	If $F=\emptyset$, we let $H_2$ be the null subgraph of $L_2$, and we let $M_2$ denote a perfect matching of~$L_2$ that contains the unique edge in $M_1\cap \partial(L)$.
	If $|F|=2$, we invoke \cref{any two adjacent edges participate in a conformal cycle and C2 is the only nonsimple theta-free mcg} and let $H_2$ denote a conformal cycle containing $F$ in $L_2$; whereas,
	if $|F|=3$, we invoke \cref{any claw centered at an isolated vertex of G-B participates in a conformal bisubdivision of theta} and let $H_2$ denote a conformal bisubdivision of $\theta$ containing $F$ in $L_2$, and in both cases, let $M_2$ be any perfect matching of $L_2-V(H_2)$.
	In all three cases, the reader may verify that $H:=H_1\cup H_2$ is a bisubdivision of $J$ in $G$, and that $M:=M_1\cup M_2$ is a perfect matching of $G-V(H)$.
\end{proof}

The reverse implication of \cref{J-freeness across a barrier where J is bicritical cubic} turns out to be more difficult to prove. We thus need a few more technical results; these, or weaker versions of them, were also used by Kothari and Murty~\cite{km16}.
The first of these is a special case of a result due to Plesn\'ik \cite[Theorem 3.4.2]{lopl86}.

\begin{lemma}\label{Plesnik - for a cc2 G: G-e-f is matchable}
	If $J$ is any $2$-connected cubic graph, then $J-e_1-e_2$ is matchable for each pair of edges $e_1$ and $e_2$.
\end{lemma}

The next result is proved in \cite{km16} for all bricks $J$, and it holds trivially if $J$ is of order two.

\begin{lemma} \label{for a nontrivial cut C of a 3-conn bicritical graph there is a pm that meets C-e in at least two edges}
	In a $3$-connected bicritical graph $J$, for any nontrivial cut $C$ and edge $e\in C$, there exists a perfect matching $M_J$ such that $|M_J\cap(C-e)|\geqslant 2$.
\end{lemma}

% The first part of the following is easily verified, whereas the second part follows from the ELP Theorem.
The following is easily verified.

\begin{proposition}\label{cubic graph with a 2-separation also has a nontrivial barrier}
	Any \con2 cubic graph $G$, that has a $2$-separation, also has a nontrivial barrier.
	Consequently, every bicritical cubic graph is \con3.
	\qed
	% Consequently, if $G$ has a nontrivial tight cut, then it has a barrier cut that is also nontrivial.
\end{proposition}

We now deduce a consequence of the above three statements.

\begin{corollary} \label{if one shore of a tightcut contains at least two branch vertices then the other shore contains at most one where J is bicritical cubic}
	Let $G$ be a \mcg{} that contains a conformal bisubdivision $H$ of a bicritical cubic graph $J$, and let $C$ be a nontrivial tight cut of $G$.
	If one shore of $C$ contains two or more branch vertices of $H$, then the other shore contains at most one branch vertex of $H$.
\end{corollary}
\begin{proof}
	Since $C$ is tight, and each perfect matching $M_H$ of $H$ extends to a perfect matching of $G$, we conclude that $|M_H\cap C|\leqslant 1$.
	In what follows, we shall invoke this observation twice.
	We let $C:=\partial(X)$, and we let $Y$~and~$\overline{Y}$ denote the sets of branch vertices of $H$ that are contained in $X$~and~$\overline{X}$, respectively.
	Our goal is to prove that either $|Y|\leqslant 1$ or $|\overline{Y}|\leqslant 1$.

	We begin by making an observation that is applicable when $Y$ and $\overline{Y}$ are both nonempty.
	Since~$J$ is \con3 by \cref{cubic graph with a 2-separation also has a nontrivial barrier}, $H$ has at least three $C$-crossing paths; by \cref{observations on C-crossing paths}~$(i)$~and~$(ii)$, each of them crosses $C$ in exactly one edge.
	If two of these $C$-crossing edges, say $e_1$~and~$e_2$, are of even parity, then \cref{Plesnik - for a cc2 G: G-e-f is matchable} implies that there exists a perfect matching~$M_H$ of $H$ that contains both $e_1$~and~$e_2$, and thus $|M_H\cap C|\geqslant 2$, contrary to our first observation.
	Hence, we conclude that at most one $C$-crossing edge is of even parity.
	If there is such an edge $e$, we let~$P_{uv}$ denote the $C$-crossing path that contains $e$, where $uv\in \partial_J(Y)$.
	Otherwise, we let $uv$ be any edge of $\partial_J(Y)$.

	Suppose to the contrary that $|Y|\geqslant 2$ and $|\overline{Y}|\geqslant 2$.
	By \cref{for a nontrivial cut C of a 3-conn bicritical graph there is a pm that meets C-e in at least two edges}, there exists a perfect matching $M_J$ of $J$ satisfying $|M_J\cap (\partial_J(Y)-uv)|\geqslant 2$;
	let $M_H$ denote the corresponding perfect matching of $H$.
	By our choice of $uv$, we infer that $|M_H\cap C|\geqslant 2$, contrary to our first observation.
	This completes the proof of \cref{if one shore of a tightcut contains at least two branch vertices then the other shore contains at most one where J is bicritical cubic}.
\end{proof}

We are now ready to prove the reverse implication of \cref{J-freeness across a barrier where J is bicritical cubic}.

% We shall now prove the reverse implication of \cref{J-freeness across a barrier where J is bicritical cubic} with the weaker hypothesis that~$J$ is a bicritical graph.

\begin{theorem}\label{if G is J-based then there is a B-fragment that is J-based where J is bicritical cubic}
	Let $G$ be a $J$-based \mcg{} that has a nonempty barrier $B$, where $J$ is a bicritical cubic graph.
	Then some \mbox{$B$-fragment} of $G$ is $J$-based.
\end{theorem}

\begin{proof}
	Let $H$ denote a conformal bisubdivision of $J$ in $G$. We shall first locate the branch vertices of~$H$, and then invoke \cref{if a tight cut shore contains at most one branch vertex then J-basedness is preserved} to finish the proof.
	\begin{statement}\label{statement:1}
		At most one branch vertex of $H$ lies in the barrier $B$.
	\end{statement} 
	\begin{proof}
		Since $J-u-v$ has a perfect matching for each pair $u,v\in V(J)$, the reader may observe that $H-u-v$ has a perfect matching, and by conformality of $H$, the graph $G-u-v$ has a perfect matching. Consequently, by \cref{when is G-u-v not matchable}, at most one of $u$ and $v$ lies in~$B$.
	\end{proof}

	% Note that the above, along with \cref{unique branch vertices of H in distinct components of G-B are nonadjacent in J} implies that if $J = \theta$, there is a $B$-fragment that contains at least one branch vertex of $H$.
	% Henceforth, we may assume that $J$ is of order four or more.

	\begin{statement}\label{statement:2}
		If $J$ has order four or more, then there exists a component of $G-B$ that contains at least two branch vertices of $H$.
	\end{statement}
	\begin{proof}
		% % We refer to the vertices of $J$ and the branch vertices of $H$ interchangeably.
		% Since $J$ has order four or more and is bicritical, each of its vertices has at least three neighbors.
		% Using \cref{statement:1}, we conclude that there exists a branch vertex of $H$, say $u$, that lies in some component $L$ of $G-B$;
		% furthermore, $u$ has at least two neighbors in $J$, say $v$ and $w$, whose corresponding branch vertices of $H$ are not in $B$.
		% If $u$ is the only branch vertex of $H$ in $L$, then by \cref{unique branch vertices of H in distinct components of G-B are nonadjacent in J}, the component of $G-B$ that contains $w$ has the desired property; otherwise, $L$ has the desired property.
		Suppose to the contrary that each $B$-fragment has at most one branch vertex of $H$.
		Then, by \cref{statement:1,unique branch vertices of H in distinct components of G-B are nonadjacent in J}, $J$ has a stable set as large as one less than its order, which is absurd.
	\end{proof}

	\begin{statement}
		There exists a component $L$ of $G-B$ that contains at least $|V(J)|-1$ branch vertices.
	\end{statement}
	\begin{proof}
		Note that by \cref{statement:1}, the statement holds for $J=\theta$.
		Otherwise $J$ has order four or more and by \ref{statement:2}, there exists a component $L$ of $G-B$ that contains at least two branch vertices of $H$.
		Observe that $\partial(L)$ is a nontrivial tight cut; since $J$ is bicritical and cubic, we invoke \cref{if one shore of a tightcut contains at least two branch vertices then the other shore contains at most one where J is bicritical cubic} to infer that $\overline{V(L)}$ contains at most one branch vertex.
	\end{proof}

	Since $\partial(L)$ is a tight cut, we invoke \cref{if a tight cut shore contains at most one branch vertex then J-basedness is preserved} to conclude that $G/\overline{V(L)}$ is $J$-based, and this completes the proof of \cref{if G is J-based then there is a B-fragment that is J-based where J is bicritical cubic}.
\end{proof}

The above, along with \cref{G is J-based if there is a B-fragment that is J-based where J has max degree 3}, proves \cref{J-freeness across a barrier where J is bicritical cubic}.
% Next, we prove \cref{J-freeness across a 2-separation where J is 3-connected cubic}.

		\subsection{Handling $2$-separations: a proof of \cref{J-freeness across a 2-separation where J is 3-connected cubic}}
			\label{sec:2sep}
			In this section, our goal is to prove \cref{J-freeness across a 2-separation where J is 3-connected cubic}.
We first prove the forward implication with weaker assumptions.

\begin{theorem} \label{if a marked S-component of G is J-based then so is G}
Let $J$ and $G$ be \mcg{}s, and let $S$ denote a $2$-separation of $G$.
If some marked $S$-component of $G$ is $J$-based, then $G$ is also $J$-based.
\end{theorem}
\begin{proof}
Let $S:=\{u,v\}$, and let $K'$ and $L'$ denote two distinct (even) components of $G-S$.
Suppose the marked $S$-component $K$ (corresponding to $K'$) is $J$-based, and let $H'$ denote a conformal bisubdivision of $J$ in $K$.
Let $e$ denote the marker edge of $K$, and let $P$ denote a conformal odd $uv$-path in $L$ provided by \cref{observations on 2-separations} $(iii)$.
Depending on whether or not $e$ belongs to $H'$, either $H:=H'+P-e$ or $H:=H'$, respectively, is a conformal bisubdivision of $J$~in~$G$.
\end{proof}

We now prove the reverse implication of \cref{J-freeness across a 2-separation where J is 3-connected cubic}, which turns out to be more difficult.

\begin{theorem} \label{if G is J-based then some marked S-component of G is J-based where J is 3-connected}
	Let $G$ be a $J$-based \mcg{} that has a $2$-separation $S$, where $J$ is a \cc3.
	Assume that if $J=\theta$ then $G-S$ has precisely two components.
	Then some marked $S$-component of $G$ is $J$-based.
\end{theorem}

\begin{proof}
Let $H$ denote a bisubdivision of $J$ that is conformal in $G$.
Since $H$ is a subdivision of a $3$-connected graph, and since $S$ is a $2$-separation of $G$, we immediately conclude that there do not exist two components of $G-S$ each of which contains a branch vertex of $H$.
This proves the following.

\begin{statement}
There exists a marked $S$-component of $G$ that contains all branch vertices of~$H$. \qed 
\end{statement} 

Since $\theta$ is the only $3$-connected cubic graph of order two, if $J\neq \theta$ then there exists a unique marked $S$-component of $G$, denoted by $L$, that contains all branch vertices of $H$.
On the other hand, if $J=\theta$ then, by our hypothesis, there are precisely two marked $S$-components; in this case, each of them may contain both branch vertices, and this happens if and only if $u$~and~$v$ are the branch vertices of $H$.
So, if $J=\theta$, we let $L$ denote the marked $S$-component that minimizes the number of $\partial_G(V(L))$-crossing paths of $H$.
Let $e$ denote the marker edge of~$L$, let $D$ denote the even cut $\partial_G(V(L))$, and let $M$ denote a perfect matching of $G-V(H)$.

\begin{statement}
	$H$ has at most one $D$-crossing path.
\end{statement}

\begin{proof}
	Suppose to the contrary that $H$ has at least two $D$-crossing paths, say $P$~and~$Q$.
	Since all branch vertices belong to $V(L)$, and since $S:=\{u,v\}$ is a $2$-separation of $G$, each of $P$~and~$Q$ contains $u$ as well as $v$.
	This immediately implies that $u$~and~$v$ are branch vertices of $H$.
	Furthermore, since $J$ is $3$-connected and cubic, we conclude that $J=\theta$.
	Consequently, $H=P\cup Q\cup R$, where $R$ is the third $uv$-path.
	Observe that if $K$ is the marked $S$-component of $G$ distinct from $L$, then the number of $\partial_G(V(K))$-crossing paths is at most one, and this contradicts our choice of $L$.
\end{proof}

We shall argue that $L$ is $J$-based.
To this end, we consider a couple of cases; in each of them, we define a subgraph $H'$ of $L$ and a matching $M'$ of $L$.

First, suppose that $H$ has no $D$-crossing paths; we let $H':=H$.
Observe that $|M\cap D|\in\{0,2\}$; if $|M\cap D|=0$, we let $M':=M\cap E(L)$, and otherwise, we let $M':=(M\cap E(L))+e$.
Now, suppose that~$H$ has precisely one $D$-crossing path, say $P_{xy}$; clearly, $|M\cap D|=0$.
We let $Q:= P_{xy}\cap (E(G)-E(L))$, let $H':=H-Q+e$ and $M':=M\cap E(L)$.
The reader may verify that in each case discussed above, $H'$ is a bisubdivision of $J$ and $M'$~is a perfect matching of $L-V(H')$.
Thus, $L$ is $J$-based, and this completes the proof of \cref{if G is J-based then some marked S-component of G is J-based where J is 3-connected}.
\end{proof}

The above, along with \cref{if a marked S-component of G is J-based then so is G}, proves \cref{J-freeness across a 2-separation where J is 3-connected cubic}.
% Next, we prove our characterization of \tf{} bricks.

		\subsection{Characterizing \tf{} bricks: a proof of \cref{K4 and Petersen are the only theta-free bricks}}
			\label{sec:bricks}
			In this section, our goal is to characterize \tf{} bricks.
Interestingly, but perhaps not surprisingly, the Petersen graph plays a crucial role in the solution to \cref{problem:theta-free}.
We start with the following relatively easy observation.

\begin{proposition}\label{Petersen is theta free but Petersen plus any edge is not}
	The Petersen graph $\mathbb{P}$ is \tf, and adding any edge to it results in a \tb{} graph.
\end{proposition}
\begin{proof}
	Clearly, one must exploit the symmetries of the Petersen graph $\mathbb{P}$.
	Since $\mathbb{P}$ has girth five and is non-hamiltonian, its only even cycles are $6$-cycles and $8$-cycles.
	If $C$ is any $8$-cycle, then $\mathbb{P}-V(C)$ is $K_2$; whereas if $C$ is any $6$-cycle, then $\mathbb{P}-V(C)$ is $K_{1,3}$.
	Thus, the conformal cycles of $\mathbb{P}$ are precisely its $8$-cycles; consequently by \cref{theta and K4 vs parity of conformal cycles}, $\mathbb{P}$ is \tf.

	Now, let $u$ and $v$ denote distinct vertices of $\mathbb{P}$.
	If $u$ and $v$ are adjacent, by \cref{any two adjacent edges participate in a conformal cycle and C2 is the only nonsimple theta-free mcg}, $\mathbb{P}+uv$ is \tb{}.
	Now suppose that $u$ and $v$ are nonadjacent; by the symmetries of $\mathbb{P}$, it suffices to consider one such pair.
	% Observe that there are two internally-disjoint odd $uv$-paths whose union is a conformal $8$-cycle $C$; thus $C+uv$ is a conformal bisubdivision of $\theta$ in $\mathbb{P}+uv$, as shown in \cref{fig:petersen plus edge}.
	Observe that there are two internally-disjoint odd $uv$-paths whose union is a conformal $8$-cycle $C$; thus $C+uv$ is a bisubdivision of $\theta$ in $\mathbb{P}+uv$, as shown in thick pink lines in \cref{fig:petersen plus edge}; it is also conformal as illustrated by the matching shown in squiggly blue line.
	% This proves the second statement.
\end{proof}

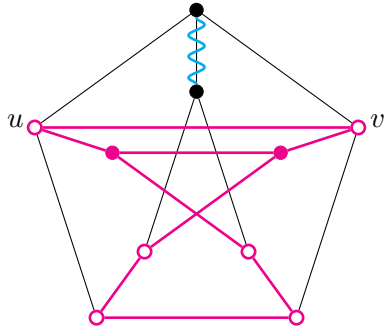
\fig{The unique simple graph obtained from the Petersen graph by adding an edge is \tb{}}{fig:petersen plus edge}{\begin{tikzpicture}[scale=1.5]

\tikzmath{
	\a = 0.78;
	\b = 1.5;
}

\node[vtx, magenta] (v0) at (18+0*72:\a) {};
\node[vtx-white,draw=magenta,thick] (u0) at (18+0*72:\b) {};

\node[vtx] (v1) at (18+1*72:\a) {};
\node[vtx] (u1) at (18+1*72:\b) {};

\node[vtx, magenta] (v2) at (18+2*72:\a) {};
\node[vtx-white,draw=magenta,thick] (u2) at (18+2*72:\b) {};

\node[vtx-white,draw=magenta,thick] (v3) at (18+3*72:\a) {};
\node[vtx-white,draw=magenta,thick] (u3) at (18+3*72:\b) {};

\node[vtx-white,draw=magenta,thick] (v4) at (18+4*72:\a) {};
\node[vtx-white,draw=magenta,thick] (u4) at (18+4*72:\b) {};

\node at (18+2*72:\b+0.18) {$u$};
\node at (18+5*72:\b+0.18) {$v$};

\draw (u0) -- (u1);
\draw (u1) -- (u2);
\draw (u2) -- (u3);
\draw (u4) -- (u0);
\draw (v4) -- (v1);
\draw (v1) -- (v3);

\draw[thick, magenta] (v0) -- (u0);
\draw (v1) -- (u1);
\draw[cyan, matching] (v1) -- (u1);
\draw[thick, magenta] (v2) -- (u2);
\draw[thick, magenta] (v3) -- (u3);
\draw[thick, magenta] (v4) -- (u4);
\draw[thick, magenta] (u3) -- (u4);
\draw[thick, magenta] (v0) -- (v2);
\draw[thick, magenta] (v2) -- (v4);
\draw[thick, magenta] (v3) -- (v0);
\draw[thick, magenta] (u2) -- (u0);

\end{tikzpicture}}

The more difficult statement to prove is that every brick that is not $K_4$ or the Petersen graph is indeed \tb{}.
To this end, we find the brick generation theorem due to Carvalho, Lucchesi and Murty \cite{clm06} useful.
In order to state their result, we need some terminology.

An edge $e$ of a matching covered graph $G$ is {\it removable} if $G-e$ is also matching covered.
In this case, if $G$ is a brick then $G-e$ need not be a brick; in particular, $G-e$ may have one or two vertices of degree two.
In order to recover a brick, at the very least, we need to get rid of the degree two vertices (if any).

Let $v_0$ denote a vertex of degree two in a matching covered graph $G$ that has two distinct neighbors, say $v_1$ and $v_2$.
Let $H$ be obtained from $G$ by contracting the two edges $v_0v_1$ and $v_0v_2$; we say that $H$ is obtained from $G$ by {\it bicontracting} the vertex $v_0$.
Observe that $H=G/X$ where $X:=\{v_0,v_1,v_2\}$ and $\partial(X)$ is a barrier cut.
The following is an immediate consequence of \cref{J-freeness across a barrier where J is bicritical cubic}.

\begin{corollary}\label{G is theta-free iff bicontraction of G is theta-free}
	Let $G$ be a matching covered graph that has a vertex $v_0$ of degree two with two distinct neighbors, and let $H$ be obtained from $G$ by bicontracting the vertex $v_0$.
	Then, $G$~is \tf{} if and only if $H$ is \tf.
	\qed
\end{corollary}

Now, let $e$ be a removable edge of a brick $G$.
As noted earlier, $G-e$ has zero, one or two vertices of degree two.
The graph $H$ obtained from $G-e$ by bicontracting each of these degree two vertices is called the {\it retract} of $G-e$.
Observe that $|V(G)|-|V(H)|\in \{0,2,4\}$; furthermore, each contraction vertex of $H$ (resulting from the bicontraction of degree two vertices), if any, has degree at least four.
We say that $e$ is a {\it thin edge} of the brick $G$ if the retract of $G-e$ is also a brick.
For instance, the bicorn shown in \cref{fig:bicorn} has a unique removable edge that is also thin.
On the other hand, the Petersen graph~$\mathbb{P}$ is devoid of thin edges since the retract of $\mathbb{P}-e$ is~$T_6$, for any edge $e$.

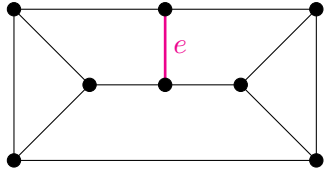
\fig{The bicorn graph and its thin edge $e$}{fig:bicorn}{\begin{tikzpicture}

\draw (-2,1)--(0,1)--(2,1)--(2,-1)--(-2,-1)--(-2,1)--(-1,0)--(-2,-1);
\draw (-1,0)--(0,0)--(1,0)--(2,1);
\draw (1,0)--(2,-1);
\draw[thick,magenta] (0,0)--(0,1);
\node[magenta] at (0.2,0.5) {$e$};
\draw (-2,1)node[vtx]{}(2,1)node[vtx]{}(-2,-1)node[vtx]{}(2,-1)node[vtx]{}(-1,0)node[vtx]{}(0,0)node[vtx]{}(1,0)node[vtx]{}(0,1)node[vtx]{};

\end{tikzpicture}}

We are now ready to state the aforementioned generation theorem of Carvalho, Lucchesi and Murty \cite{clm06}, a stronger version of a conjecture of Lov\'asz proved earlier by the same authors~\cite{clm02, clm02a}.

\begin{theorem}{\sc [Thin Edge Theorem]}\label{CLM thin edge theorem}
\newline
	Every brick distinct from $K_4$, $\overline{C_6}$ and the Petersen graph $\mathbb{P}$ has a thin edge.
\end{theorem}

In their work \cite{clm06}, Carvalho, Lucchesi and Murty present the above theorem from the generation viewpoint, stating that every brick may be generated from at least one of $K_4, \overline{C_6}$ and $\mathbb{P}$ by means of a sequence of `expansion operations'.
They defined four expansion operations --- each of which comprises a sequence of zero, one or two `bisplittings' (the inverse of special types of bicontractions) followed by adding an edge.
In particular, if $e$ is a thin edge of a brick $G$, and $H$ is the retract of $G-e$, then~$G$ may be obtained from $H$ using one of the four expansion operations.
We do not require the precise definitions of these operations to prove our characterization of \tf{} bricks restated below; the interested reader may refer to \cite{clm06,lumu24}.

\tfBricks*

\begin{proof}
	Clearly, $K_4$ is \tf{} and by \cref{Petersen is theta free but Petersen plus any edge is not}, $\mathbb{P}$ is \tf.
	Conversely, let $G$ be a \tf{} brick.
	We proceed by induction on the size of $G$.
	First suppose that $G$ has a thin edge, say $e$, and let $H$ denote the retract of $G-e$.
	Clearly, $G-e$ is \tf, and by \cref{G is theta-free iff bicontraction of G is theta-free}, the brick $H$ is also \tf;
	by the induction hypothesis, $H$ is either $K_4$ or $\mathbb{P}$.
	As noted earlier, any contraction vertex of $H$ that results from the bicontraction of degree two vertices of $G-e$ has degree at least four.
	Since $H$ is cubic, we conclude that $H$ is the same as $G-e$; whence, $G$ may be obtained from $H$ by adding an edge.
	Depending on whether $H$ is $K_4$ or $\petersen$, we invoke \cref{any two adjacent edges participate in a conformal cycle and C2 is the only nonsimple theta-free mcg} or \cref{Petersen is theta free but Petersen plus any edge is not}, respectively, to infer that $G$ is \tb{}; a contradiction.
	By the Thin Edge Theorem (\ref{CLM thin edge theorem}) and the fact that $\overline{C_6}$ is \tb{}, the desired conclusion holds.
\end{proof}

This concludes the proof of our Main Theorem (\ref{theta-free characterization - inductive version}).

	\section{Consequences of the Main Theorem}
	\label{sec:consequences of main theorem}
		In this section, we prove all of the consequences of the Main Theorem discussed in \cref{subsec:consequences}.

		\subsection{Sparsity and bicriticality}
			\label{subsec:sparsity and bicriticality}
			Recall the splicing operation from \cref{sec:characterizing theta-free graphs}.
The following is easily verified; see \cite[Proposition 5.14]{lumu24}.

\begin{proposition}
	\label{splicing two bicritical graphs yields a biciritical graph}
	Splicing biciritical graphs yields a biciritical graph.
	\qed
\end{proposition}

Also, recall our discussion from \cref{subsubsec:ELP cuts} that compares marked $S$-components of a \mcg{} $G$ with its $C$-contractions, where $C$ is any $2$-separation cut arising from a \mbox{$2$-separation $S$}.
This, along with the above proposition, and the fact that a graph is bicritical if and only if its underlying simple graph is bicritical, immediately implies the reverse implication of the following; the forward implication may be verified easily.

\begin{proposition}
	\label{bicritical iff k2sum is bicritical}
	Let $G$ be any $K_2$-sum of graphs $G_1$ and $G_2$, neither of which is isomorphic to~$K_2$.
	Then, $G$ is bicritical if and only if both $G_1$ and $G_2$ are bicritical.
	\qed
\end{proposition}

% We now show that the converse of the above statement also holds.
% The following observation on the \mbox{$K_2$-sum} operation is crucial for characterizing biciritical \tf{} graphs.

% \unsure{is there a known result that gives us this with a shorter (or no) proof}
% \begin{proof}
% 	Let $G$ be the $K_2$-sum of $G_1$ and $G_2$ at the edges $x_1y_1$ and $x_2y_2$, respectively.
% 	Let $x$ and $y$ denote the vertices of $G$ obtained by identifying $x_i$'s and $y_i$'s, respectively, for $i \in \{1,2\}$.

% 	First, suppose that $G$ is bicritical. Let $u,v \in V(G_1)$ and let $M$ denote a perfect matching of $G-u-v$. Since $|V(G_1)-u-v|$ is even, $|M \cap \partial_{G-u-v}(V(G_1)-u-v)|$ is even; in particular, it is either zero or two since $\{x,y\}$ is a vertex cut. In the first case, let $M_1:=M \cap E(G_1)$ and in the other case, let $M_1:=M \cap E(G_1) + x_1y_1$. Note that in either case, $M_1$ is a perfect matching of $G_1-x-y$ and hence, $G_1$ is bicritical; by the symmetry, $G_2$ is also bicritical.

	% Now, suppose $G_1$ and $G_2$ are bicritical. Let $u,v$ be any distinct vertices of $V(G)$. Adjust notation so that $u \in V(G_1)$. If $v \in G_1$ as well, then let $M_1$ be a perfect matching of $G_1-u-v$ and let $M_2$ be a perfect matching of $G_2-x_2-y_2$. Otherwise, if $v \in V(G_2)$, adjust notation so that $u \neq x$, and let $M_1$ be a perfect matching of $G_1-u-x_1$ and let $M_2$ be a perfect matching of $G_2-v-y_2$. Observe that in each case, $M:=M_1 \cup M_2$ is a perfect matching of $G-u-v$ and hence, $G$ is bicritical.
% \end{proof}

The above, along with \cref{theta-free characterization - inductive version}, immediately implies the following characterization of \tf{} bicritical graphs.

\begin{restate}{\ref*{characterization of theta-free-bicritical graphs}}
	{\sc[\tf{} Bicritical Graphs]} \newline
	A \tf{} graph $G$ is bicritical if and only if $G \in \{K_2\} \cup \mathcal{T}$.
	\qed
\end{restate}

% not using the below command in order to add a qed
% \tfBicritical*

We now switch our attention to proving upper bounds on the size of any \tf{} graph.
We start by discussing a couple of auxiliary results that are easily observed.
The aforementioned discussion on marked $S$-components versus $C$-contractions, coupled with Lov\'asz's Unique Tight Cut Decomposition Theorem (\ref{unique tcd theorem}), implies the following.

\begin{proposition}\label{b(G) across 2-separations}
	For any \mcg{} $G$ that has a $2$-separation $S$, the following holds:
	\[
		b(G)=\sum_{L\in\mathcal{L}} b(L)
	\]
	where $\mathcal{L}$ denotes the collection of its marked $S$-components. \qed
\end{proposition}

Observe that the bricks and braces of a matching covered graph $G$, that has a nontrivial barrier~$B$, is the union of the bricks and braces of its \mbox{$B$-fragments} (of order at least four) as well as the braces of its $B$-core.
This, along with Lov\'asz's Unique Tight Cut Decomposition Theorem~(\ref{unique tcd theorem}), proves the following.

\begin{proposition}\label{b(G) across a barrier}
	For any \mcg{} $G$ that has a nonempty barrier $B$, the following holds:
	\[
		b(G)=\sum_{i=1}^{|B|} b(G_i)
	\]
	where $G_1,G_2,\dots,G_{|B|}$ denote its \mbox{$B$-fragments}. \qed
\end{proposition}

% Recall that if $S:=\{u,v\}$ is a $2$-separation of a \mcg{}~$G$, and $L$ is a component of $G-S$, the marked $S$-component $L'$ corresponding to $L$, and the tight cut contraction $G/(\overline{V(L)\cup \{u\}}\rightarrow v)$ are isomorphic except for the multiplicities of the edges joining $u$~and~$v$.

We are now ready to prove \cref{edge bound on theta-free graphs}.

\edgeBoundTheorem*

\begin{proof}
	We shall prove statements~$(i)$~and~$(ii)$ by induction on the order.
	Let $G$ be any \tf{} graph;
	we invoke \cref{theta-free characterization - inductive version}.
	If statement $(i)$ of \cref{theta-free characterization - inductive version} holds, the desired conclusions hold.
	If statement $(ii)$ of \cref{theta-free characterization - inductive version} holds, let $B$ denote a nontrivial barrier, and if statement $(iii)$ holds, let $S$ denote a $2$-separation.
	In the former case, we let~$G_1,G_2,\dots,G_{|B|}$ denote the $B$-fragments, whereas in the latter case, we let~$G_1$~and~$G_2$ denote the marked $S$-components;
	in either case, by the induction hypothesis, each $G_i$ satisfies $m_i \leqslant \frac{3n_i}{2}+b_i-1$ and $b_i \leqslant \frac{n_i}{2}-1$, where $n_i,m_i$ and $b_i$ denote its order, size and the number of bricks.

	If $G$ has a nontrivial barrier $B$, using \cref{b(G) across a barrier}, the fact that the edge sets of the \mbox{$B$-fragments} comprise a partition of $E(G)$, and the fact that the sum of the orders of the \mbox{$B$-fragments} equals $n$:
	\[m=\sum_{i=1}^{|B|} m_i\leqslant \sum_{i=1}^{|B|} \left( \frac{3n_i}{2}+b_i-1\right)=\frac{3n}{2}+b-{|B|}< \frac{3n}{2}+b-1, ~\,\text{and}\]
	\[b=\sum_{i=1}^{|B|}b_i \leqslant \sum_{i=1}^{|B|} \left(\frac{n_i}{2}-1\right)=\frac{n}{2}-|B|<\frac{n}{2}-1\]
	Therefore, each bound holds strictly.
	By \cref{when is G-u-v not matchable}, $G$ is not bicritical; consequently, by \cref{characterization of theta-free-bicritical graphs}, $G \notin \mathcal{T}$.
	Thus, the desired conclusion holds.

	On the other hand, if $G$ has a $2$-separation $S$, observe that $n = n_1 + n_2 - 2$ and that $m = m_1 + m_2 - 2$.
	Therefore, by \cref{b(G) across 2-separations}:
	\[m = m_1 + m_2 -2 \leqslant \left(\frac{3n_1}{2}+b_1-1\right)+\left(\frac{3n_2}{2}+b_2-1\right)-2=\frac{3(n+2)}{2}+b-4=\frac{3n}{2}+b-1, ~\,\text{and}\]
	\[b=b_1+b_2\leqslant \left(\frac{n_1}{2}-1\right)+\left(\frac{n_2}{2}-1\right)=\frac{n+2}{2}-2=\frac{n}{2}-1\]
	Thus, both desired inequalities hold for $G$;
	furthermore, either of them holds with equality if and only if the corresponding inequality holds with equality for each of $G_1$ and $G_2$.
	By the induction hypothesis, and the fact that $G$ is a $K_2$-sum of~$G_1$~and~$G_2$, the desired conclusions hold.
\end{proof}

We conclude this section by defining a family $\mathcal{K}$ (that is strictly contained in $\mathcal{T}_0$) that demonstrates \cref{for any k there is a tf graph with a vertex of degree k}.
The brick $K_4$ belongs to $\mathcal{K}$, and
any $K_2$-sum of $K_4$ and a member of $\mathcal{K}$, at an edge that is incident with a vertex of maximum degree,
also belongs to $\mathcal{K}$.
One may easily verfiy that the maximum degree of any member of~$\mathcal{K}$ is $\frac{n}{2}+1$.

		\subsection{\tf{} near-bipartite graphs}
			\label{subsec:tf near-bipartite}
			Recall that the number of bricks of a \mcg{} $G$ is denoted by $b(G)$.
We start with the following basic fact about near-bipartite graphs \cite[Exercise 8.3.3]{lumu24}.

\begin{proposition}
	\label{a near-brick is near-bipartite if and only if its unique brick is near-bipartite}
	A \mcg{} $G$ is near-bipartite if and only if $b(G)=1$ and its unique brick is near-bipartite.
\end{proposition}

Recall that a \mcg{} is cycle-extendable if each of its even cycles is conformal.
A result of Carvalho and Little \cite[Corollary 2.4]{cali08} implies the following.

\begin{proposition}
	\label{b+p <= 1 for cycle-extendable graphs}
	Every nonbipartite cycle-extendable graph $G$ satisfies $b(G) = 1$ and its unique brick is not the Petersen graph $\petersen$.
\end{proposition}

We are now ready to prove our characterization of cycle-extendable \tf{} graphs.

\tfNearBipartite*

\begin{proof}
	Let $G$ be a \tf{} graph.
	We may assume that $G$ is nonbipartite; otherwise, the result follows immediately from \cref{the only bipartite theta-free mcgs are k2 and even cycles}.
	We shall show: $(i) \implies (iv) \implies (iii) \implies (ii) \implies (i)$.

	First suppose that $(i)$ holds and let $\{P,Q\}$ be a removable double ear such that $G-P-Q$ is bipartite.
	By \cref{the only bipartite theta-free mcgs are k2 and even cycles}, $G-P-Q$ is an even cycle.
	The reader may verify using the definition of removable double ears that $G$ is indeed a bisubdivision of $K_4$, whence $(iv)$ holds.
	The reader may easily verify that $(iv)$ implies $(iii)$.
	Now suppose that $(iii)$ holds and note that by \cref{each brick of a tf graph is either K4 or the Petersen graph,b+p <= 1 for cycle-extendable graphs}, $b(G)=1$ and the unique brick of $G$ is $K_4$; whence, $(ii)$ holds.
	Finally, if $(ii)$ holds, by \cref{a near-brick is near-bipartite if and only if its unique brick is near-bipartite}, $G$ is near-bipartite and $(i)$ holds.
\end{proof}

		\subsection{Pfaffian \tf{} graphs}
			\label{subsec:tf Pfaffian}
			We start with some terminology required to define the term `$\sminor{}$'.
A cut $C$ of a \mcg{} is \emph{separating} if both $C$-contractions are matching covered.
For instance, tight cuts are special types of separating cuts.
However, not all separating cuts are tight; consider the cut $C$ shown in \cref{fig:c6bar}.

A \emph{removable class} of a \mcg{} is any minimal nonempty subset of its edges whose deletion yields a \mcg.
For example, the singleton set containing any edge is a removable class of $K_{3,3}$; however, $K_4$ has no such singletons but each perfect matching is a removable class.
Lucchesi and Murty \cite[Corollary 8.18]{lumu24} proved that any removable class of a \mcg{} has cardinality one or two; this is closely related to the notion of removable ears defined in \cref{subsec:tf near-bipartite intro}.
If a removable class $R$ is a singleton, its only member is a removable edge;
otherwise, $R$ is a \emph{removable doubleton}.

A \mcg{} $G$ contains a (matching covered) graph $J$ as a \emph{separation-deletion minor}, or simply as an \emph{\sminor}, if $J$ may be obtained from $G$ by a sequence of operations, each of which is either shrinking a shore of a separating cut or deleting a removable class.
This containment notion is weaker than that of the conformal minor;
for instance, each of the graphs shown in \cref{fig:barrier} contains $K_{3,3}$ as an \sminor{}, but all of them are $K_{3,3}$-free.
The following result \cite[Lemma 21.1]{lumu24} shows that the \sminor{} containment notion preserves the Pfaffian property.

\begin{theorem}\label{every S-minor of a Pfaffian graph is Pfaffian}
	A \mcg{} is Pfaffian if and only if each of its \sminor{}s is Pfaffian.
\end{theorem}

Due to the above, one approach to characterize Pfaffian graphs is to list down all the \emph{$S$-minimally non-Pfaffian graphs}, that is, those non-Pfaffian \mcg{}s, each of whose proper \sminor{}s is Pfaffian.
This is precisely what happened over the last few decades for various graph classes.
Before we discuss these results, we state the following \cite[Theorem 20.10]{lumu24}, which implies that all $S$-minimally non-Pfaffian graphs are either bricks or braces.

\begin{theorem}\label{G is Pfaffian iff tight cut contractions of G are Pfaffian}
	A matching covered graph $G$ is Pfaffian if and only if both of its $C$-contractions are Pfaffian, where $C$ is any tight cut of $G$.
	Consequently, $G$ is Pfaffian if and only if each of its bricks and braces is Pfaffian.
\end{theorem}

The forward implication of the above is a special case of \cref{every S-minor of a Pfaffian graph is Pfaffian}.
As mentioned in \cref{subsubsec:tf pfaffian intro}, Little \cite{litt75} showed that $K_{3,3}$ is the only minimally non-Pfaffian bipartite graph; recall \cref{characterization of bipartite Pfaffian graphs}.
Fischer and Little proved the following characterization of Pfaffian near-bipartite graphs.

\begin{theorem}\label{near-bipartite Pfaffian}
	For a near-bipartite graph $G$, the following are equivalent:
	\begin{enumerate}[label=(\roman*)]
		\item $G$ is Pfaffian,
		\item $G$ does not contain any of $K_{3,3},\,\Gamma_1$ and $\Gamma_2$ as an \sminor, and
		\item $G$ is \mbox{$J$-free} for each $J \in \{K_{3,3},~\Gamma_1,~\Gamma_2,~\Gamma_{0,1},~\Gamma_{0,2},~\Gamma_{1,2}\}$.
	\end{enumerate}
\end{theorem}	

The graphs $\Gamma_{0,1}$ and $\Gamma_{0,2}$ are the ones shown in \cref{fig:barriertb1,fig:k33 spliced with 2 k4s on either side}, respectively;
for the rest, we refer the reader to \cite{fili01}.
The original statement \cite[Theorem 1.3]{fili01} lists seven graphs;
however, one of them, namely $\Gamma_{1,1}$, is $K_{3,3}$-based and thus redundant.

Recall that separating cuts are not necessarily tight; this inspires the following definition.
A \mcg{} $G$ is \emph{solid} if each of its separating cuts is also tight.
It is worth noting that bipartite \mcg{}s comprise a proper subset of the solid ones.
We refer the interested reader to \cite[Chapter 7]{lumu24} for further discussion on solid graphs and their importance.
Carvalho, Lucchesi and Murty \cite{clm12} generalized Little's Theorem (\ref{characterization of bipartite Pfaffian graphs}) as follows.

\begin{theorem}\label{solid Pfaffian}
	For a solid \mcg{} $G$, the following are equivalent:
	\begin{enumerate}[label=(\roman*)]
		\item $G$ is Pfaffian,
		\item $G$ does not contain $K_{3,3}$ as an \sminor, and
		\item $G$ is $K_{3,3}$-free, and $J$-free for each graph $J$ shown in \cref{fig:barrier}.
	\end{enumerate}
\end{theorem}

We remark that our characterization of Pfaffian \tf{} graphs (\cref{theta-free Pfaffian characterization}) is in the same spirit as all of the aforementioned results (\cref{characterization of bipartite Pfaffian graphs,near-bipartite Pfaffian,solid Pfaffian}) that characterize various subclasses of Pfaffian graphs.
We now discuss a few results that we shall find useful in proving it.
The ELP Theorem (\ref{ELP theorem - bricks version}), along with \cref{cubic graph with a 2-separation also has a nontrivial barrier,splicing two bicritical graphs yields a biciritical graph}, immediately implies the following.

\begin{proposition}\label{splicing two cubic bricks yeilds a cubic brick}
Any splicing of two cubic bricks yields a cubic brick. \qed
\end{proposition}

\fig{A \mbox{$K_4$-decoration} of $K_{3,3}$}{fig:k33 spliced with 2 k4s on either side}{\begin{tikzpicture}[scale=1.5]

	\tikzmath{
		\height = 0.6;
		\halfsidelength = sin(30)*\height;
		\sidelength = 2*\halfsidelength;
		\x = 1.25*2;
		\y = -1.125;
		\offset = 0.2;
	}

	% vertices
	\node[vtx] (a22) at (\x,\y+2*\height/3-\offset) {};
	\node[vtx] (a21) at (\x-\halfsidelength, \y-\height/3-\offset) {};
	\node[vtx] (a23) at (\x+\halfsidelength, \y-\height/3-\offset) {};

	\node[vtx] (b22) at (\x, -2*\height/3+\offset) {};
	\node[vtx] (b21) at (\x-\halfsidelength, \height/3+\offset) {};
	\node[vtx] (b23) at (\x+\halfsidelength, \height/3+\offset) {};

	\foreach \i in {1,3}{
		\node[vtx] (b\i) at (1.25*\i,0) {};
		\node[vtx] (a\i) at (1.25*\i,\y) {};
	}

	% edges
	\draw (a21) -- (a22) -- (a23) -- (a21);
	\draw (b21) -- (b22) -- (b23) -- (b21);

	\draw (a21) -- (b1);
	\draw (a22) -- (b22);
	\draw (a23) -- (b3);

	\draw (b21) -- (a1);
	\draw (b23) -- (a3);

	\draw (a1) -- (b3);
	\draw (a3) -- (b1);

	\draw (a1) -- (b1);
	\draw (a3) -- (b3);

\end{tikzpicture}}

Let $X$ denote any (possibly empty) subset of the vertices of a \con2 cubic graph $J$.
The graph obtained from $J$, by splicing at each vertex in $X$ with a copy of $K_4$, is called a \emph{\mbox{$K_4$-decoration}} of $J$.
An example is shown in \cref{fig:k33 spliced with 2 k4s on either side}; observe that this graph is a brick.
Using this, and \cref{K4 and Petersen are the only theta-free bricks,splicing two cubic bricks yeilds a cubic brick}, we conclude that for any \tf{} $K_4$-decoration of $K_{3,3}$, the set~$X$ is a subset of one of the color class of $K_{3,3}$.
This, along with the facts illustrated in \cref{fig:barrier}, implies the second part of the following;
the first part follows immediately from \cref{splicing two cubic bricks yeilds a cubic brick,K4 and Petersen are the only theta-free bricks}.

\begin{proposition}\label{tf k4-decorations of k33 and petersen}
	The Petersen graph $\petersen$ is the only \tf{} $K_4$-decoration of itself, whereas $\tftwelve$ is the only \tf{} $K_4$-decoration of $K_{3,3}$.
	\qed
\end{proposition}

The following unpublished result implies that every class with finitely many forbidden \sminor{}s, that are all \con2 cubic graphs, has finitely many forbidden conformal minors that are also \con2 and cubic; see  \cite[Lemma 1.18]{klls26}.

\begin{lemma}\label{k4-decoration lemma}
	A matching covered graph $G$ contains a \cc2 $J$ as an \sminor{} if and only if it contains some \mbox{$K_4$-decoration} of $J$ as a conformal minor.
\end{lemma}

We are now ready to prove our characterization of Pfaffian \tf{} graphs.

\tfPfaffian*

\begin{proof}
	Let $G$ be a \tf{} graph.
	Using \cref{every S-minor of a Pfaffian graph is Pfaffian}, and the fact that $K_{3,3}$ and $\petersen$ are non-Pfaffian, we deduce that $(i)$ implies $(ii)$.
	Now suppose that $(ii)$ holds.
	Note that each brick and brace of $G$ is also an \sminor{};
	hence, $\petersen$ is not a brick of $G$.
	Furthermore, by invoking the equivalence of $(i)$ and $(ii)$ of \cref{characterization of bipartite Pfaffian graphs} to each brace $H$ of $G$, and by transitivity of \sminor{} containment, we conclude that $H$ is Pfaffian.
	Thus, $(ii)$ implies $(iv)$.
	Now, suppose that $(iv)$ holds.
	Then, by \cref{each brick of a tf graph is either K4 or the Petersen graph}, each brick of $G$ is the Pfaffian graph $K_4$, whence by \cref{G is Pfaffian iff tight cut contractions of G are Pfaffian}, $G$ is Pfaffian;
	thus, $(iv)$ implies $(i)$.
	Finally, the fact that $(ii)$ and $(iii)$ are equivalent follows immediately from \cref{k4-decoration lemma,tf k4-decorations of k33 and petersen}.
\end{proof}

		\subsection{Cubic \tf{} graphs}
			\label{subsec:tf cubic}
			In this section, we shall view cycles as sets of edges.
Observe that a conformal cycle $Q$ of a \mcg{} meets each tight cut in zero or two edges; in the latter case, these two edges belong to different perfect matchings of $Q$.
Using this, the forward implication of the following is easy to see; the reverse implication is straightforward.

\begin{proposition}\label{conformal cycles across tight cuts}
	A cycle $Q$ of a matching covered graph $G$ is conformal if and only if its restriction to each \mbox{$C$-contraction} is either a conformal cycle or empty, where $C$ is any tight cut of~$G$. \qed
\end{proposition}

The following is an immediate consequence.
For instance, consider the cycle $Q$, shown in thick orange lines, in \cref{fig:barriertf}.

\begin{corollary}
	\label{conformal cycles along a barrier}
	A cycle $Q$ of a matching covered graph $G$ is conformal if and only if, for each $B$-fragment $J$, the restriction of $Q$ to $J$ is either empty or a conformal cycle of $J$, where $B$ is any barrier of $G$. \qed
\end{corollary}

We are now ready to prove \cref{a cubic graph is theta free iff each conformal cycle is of length 0 mod 4}.

\tfCubic*

\begin{proof}
	The forward implication is simply the first part of \cref{theta and K4 vs parity of conformal cycles}.
	Conversely, let $G$ be a \cc2 that is \tf{}.
	We proceed by induction on the order.
	If statement~$(i)$ of \cref{theta-free characterization - inductive version} holds, then $G$ is either $K_4$ or $\petersen$; in the former case, each conformal cycle is of length four, whereas in the latter case, each conformal cycle is of length eight.
	Otherwise, by \cref{cubic graph with a 2-separation also has a nontrivial barrier}, statement~$(ii)$ of \cref{theta-free characterization - inductive version} holds;
	that is, $G$ has a nontrivial barrier $B$ and each of its \mbox{$B$-fragment} is \tf{}.
	A straightforward counting argument shows that each $B$-fragment is cubic;
	using the fact that vertex-connectivity equals edge-connectivity, we deduce that each of them is \con2.
	Now, let $Q$ be any conformal cycle of $G$.
	Since the edge sets of the $B$-fragments comprise a partition of $E(G)$, each edge of $Q$ belongs to precisely one \mbox{$B$-fragment}.
	Invoking \cref{conformal cycles along a barrier} and the induction hypothesis, we arrive at the desired conclusion.
	% Let $Q$ be any conformal cycle of $G$; let $M$ be a perfect matching of $G-V(Q)$.
	% For $i \in \{1,2,\dots,|B|\}$, let $Q_i:=Q \cap E(G_i)$.
	% Since $Q$ is a cycle, $|Q \cap \partial_G(J_i)|$ is even and by \cref{every tight cut of a cubic graph is a 3-cut}, $\partial_G(J_i)$ is a $3$-cut.
	% Therefore, $|Q \cap \partial_G(J_i)|$ is either zero or two and consequently, $C_i$ is either empty or a cycle for each $i \in \{1,2,\dots,|B|\}$.
	% If $C_i$ is a cycle, note that $M \cap E(G_i)$ is a perfect matching of $G_i$ and hence, $C_i$ is a conformal cycle of $G_i$.
	% By the induction hypothesis, $C_i$ is of length $0 \pmod{4}$.
	% Observe that $C$ is exactly the disjoint union of $C_i$'s for $i \in \{1,2,\dots,|B|\}$ and hence $C_i$ is of length $0 \pmod{4}$.
\end{proof}

		\subsection{$3$-edge-colorability of \tf{} cubic graphs}
			\label{subsec:tf cubic 3ec}
			% For a cubic graph, its vertex-connectivity is the same as its edge-connectivity.
% Let $G$ be any \con2 cubic graph.
% If $G$ has a \mbox{$2$-cut} $C$ then $G-C$ has precisely two components;
% the (cubic) graph obtained from either component, by adding an edge joining its vertices of degree two, is a \emph{marked $C$-component} of $G$.
% The following is easily verified.

% \begin{proposition}
% 	A \cc2 $G$ is $3$-edge-colorable if and only if both of its marked $C$-components are $3$-edge-colorable, where $C$ is any $2$-cut of $G$. \qed
% \end{proposition}

% The above implies that, for \cref{dp:cubic-3ec}, one may restrict attention to \cc3s;
% we now go a step further.
Using the fact that every perfect matching meets each odd cut in at least one edge, one may observe the forward implication of the following; the reverse is straightforward.

\begin{proposition}
	A \cc2 $G$ is $3$-edge-colorable if and only if both of its \mbox{$C$-contractions} are $3$-edge-colorable, where $C$ is any $3$-cut of $G$. \qed
\end{proposition}

Combining the above, with the facts that $(i)$ each tight cut in a \cc2 is a $3$-cut \cite[Theorem 5.8]{lumu24} and $(ii)$ bipartite cubic graphs are $3$-edge-colorable, we infer the following.

\begin{corollary}
	\label{cc2 is 3-e-c iff each of its bricks is}
	A \cc2 is $3$-edge-colorable if and only if each of its (cubic) bricks is $3$-edge-colorable. \qed
\end{corollary}

We are now ready to prove our characterization of \mbox{$3$-edge-colorable} \tf{} cubic graphs.

\tfCubicThreeEdgeColorability*

\begin{proof}
	Let $G$ be a \tf{} cubic graph.
	Note that $K_4$ is $3$-edge-colorable but the Petersen graph is not.
	Using this, \cref{cc2 is 3-e-c iff each of its bricks is,each brick of a tf graph is either K4 or the Petersen graph} we conclude that $(i), (ii)$ and $(iii)$ are equivalent.
	The equivalence of $(iii)$ and $(iv)$ follows from the fact that the Petersen graph contains itself as a conformal minor and \cref{for a cubic brick J: G is J-free iff each brick of G is J-free}.
	The equivalence of $(iv)$ and $(v)$ follows from \cref{k4-decoration lemma,tf k4-decorations of k33 and petersen}.
\end{proof}

% old version below

% Finally, using the above, \cref{each brick of a tf graph is either K4 or the Petersen graph} and the facts that $K_4$ is $3$-edge-colorable but the Petersen graph is not, we deduce the following.

% \begin{restate}{\ref*{3-e-c cubic theta free graphs}}
% 	{\sc[$3$-edge-colorability of \tf{} Cubic Graphs]}\newline
% 	For a \tf{} cubic graph $G$, the following are equivalent:
% 	\begin{enumerate}[label=(\roman*)]
% 		\item $G$ is $3$-edge-colorable,
% 		\item each brick of $G$ is isomorphic to $K_4$,
% 		\item the Petersen graph is not a brick of $G$,
% 		\item $G$ does not contain the Petersen graph as a conformal minor, that is, $G$ is \mbox{$\petersen$-free}, and
% 		\item $G$ does not contain the Petersen graph as an \sminor{}. \qed
% 	\end{enumerate}
% \end{restate}

	% \bibliographystyle{plain}
	% \bibliography{clm}

\begin{thebibliography}{10}

\bibitem{bomu08}
J.~A. Bondy and U.~S.~R. Murty.
\newblock {\em Graph Theory}.
\newblock Springer, 2008.

\bibitem{cali08}
M.~H. Carvalho and C.~H.~C. Little.
\newblock Vector spaces and the {P}etersen graph.
\newblock {\em The Electronic Journal of Combinatorics}, 15(1):\#R9, Jan. 2008.

\bibitem{clm02}
M.~H. Carvalho, C.~L. Lucchesi, and U.~S.~R. Murty.
\newblock On a {C}onjecture of {L}ov{\'a}sz {C}oncerning {B}ricks. {I}. {T}he
  {C}haracteristic of a {M}atching {C}overed {G}raph.
\newblock {\em J.~Combin.~Theory Ser.~B}, 85:94--136, 2002.

\bibitem{clm02a}
M.~H. Carvalho, C.~L. Lucchesi, and U.~S.~R. Murty.
\newblock On a {C}onjecture of {L}ov{\'a}sz {C}oncerning {B}ricks. {I}{I}.
  {B}ricks of {F}inite {C}haracteristic.
\newblock {\em J.~Combin.~Theory Ser.~B}, 85:137--180, 2002.

\bibitem{clm05a}
M.~H. Carvalho, C.~L. Lucchesi, and U.~S.~R. Murty.
\newblock On the number of dissimilar {P}faffian orientations of graphs.
\newblock {\em {RAIRO} - Inf.~Theor.~Appl.}, 39:93--113, 2005.

\bibitem{clm06}
M.~H. Carvalho, C.~L. Lucchesi, and U.~S.~R. Murty.
\newblock How to build a brick.
\newblock {\em Discrete Math.}, 306:2383--2410, 2006.

\bibitem{clm12}
M.~H. Carvalho, C.~L. Lucchesi, and U.~S.~R. Murty.
\newblock A generalization of {L}ittle's {T}heorem on {P}faffian orientations.
\newblock {\em J.~Combin.~Theory Ser.~B}, 102:1241--1266, 2012.

\bibitem{clm18}
M.~H. Carvalho, C.~L. Lucchesi, and U.~S.~R. Murty.
\newblock On tight cuts in matching covered graphs.
\newblock {\em Journal of Combinatorics}, 2018.

\bibitem{cfllz21}
G.~Chen, X.~Feng, F.~Lu, C.~L. Lucchesi, and L. Zhang.
\newblock Laminar tight cuts in matching covered graphs.
\newblock {\em J.~Combin.~Theory Ser.~B}, 150:177--194, 2021.

\bibitem{dpdk25}
A.~Y. Dalwadi, K.~R.~S. Pause, A.~A. Diwan, and N.~Kothari.
\newblock Planar cycle-extendable graphs.
\newblock {\em Discrete Mathematics \& Theoretical Computer Science}, vol.
  27:2, May 2025.

\bibitem{elp82}
J.~Edmonds, L.~Lov\'asz, and W.~R. Pulleyblank.
\newblock Brick decomposition and the matching rank of graphs.
\newblock {\em Combinatorica}, 2:247--274, 1982.

\bibitem{fili01}
I.~Fischer and C.~H.~C. Little.
\newblock A {C}haracterisation of {P}faffian {N}ear {B}ipartite {G}raphs.
\newblock {\em J.~Combin.~Theory Ser.~B}, 82:175--222, 2001.

\bibitem{fjls03}
M.~Funk, B.~Jackson, D.~Labbate, and J.~Sheehan.
\newblock Det-{E}xtremal {C}ubic {B}ipartite {G}raphs.
\newblock {\em Journal of Graph Theory}, 44(1):50--64, 2003.

\bibitem{holy81}
I.~Holyer.
\newblock The {N}{P}-{C}ompleteness of {E}dge-{C}oloring.
\newblock {\em SIAM Journal on Computing}, 10(4):718--720, 1981.

\bibitem{kast63}
P.~W. Kasteleyn.
\newblock Dimer {S}tatistics and {P}hase {T}ransitions.
\newblock {\em J. Math. Phys.}, 4:287--293, 1963.

\bibitem{koth16}
N.~Kothari.
\newblock {\em Brick Generation and Conformal Subgraphs}.
\newblock PhD thesis, University of Waterloo, 2016.
\newblock Available at \url{https://uwspace.uwaterloo.ca/handle/10012/10376}.

\bibitem{klls26}
N.~Kothari, O.~Lee, C.~L. Lucchesi, and C.~N. Silva.
\newblock Cubic graphs, ${S}$-minors and conformal minors.
\newblock Preprint available at \url{https://arxiv.org/abs/2606.04173}, 2026.

\bibitem{km16}
N.~Kothari and U.~S.~R. Murty.
\newblock ${K}_{4}$-free and $\overline{C_6}$-free {P}lanar {M}atching
  {C}overed {G}raphs.
\newblock {\em J.~Graph Theory}, 82(1):5--32, 2016.

\bibitem{gale83}
D.~Leven and Z.~Galil.
\newblock N{P} {C}ompleteness of {F}inding the {C}hromatic {I}ndex of {R}egular
  {G}raphs.
\newblock {\em Journal of Algorithms}, 4(1):35--44, 1983.

\bibitem{litt75}
C.~H.~C. Little.
\newblock A {C}haracterization of {C}onvertible $(0,1)$-{M}atrices.
\newblock {\em J.~Combin.~Theory Ser.~B}, 18:187--208, 1975.

\bibitem{lire91}
C.~H.~C. Little and F.~Rendl.
\newblock Operations preserving the {P}faffian property of a graph.
\newblock {\em J. Austral. Math. Soc. (Series A)}, 50:248--275, 1991.

\bibitem{lova83}
L.~Lov\'asz.
\newblock Ear-decompositions of matching-covered graphs.
\newblock {\em Combinatorica}, 3:105--117, 1983.

\bibitem{lova87}
L.~Lov\'asz.
\newblock Matching {S}tructure and the {M}atching {L}attice.
\newblock {\em J.~Combin.~Theory Ser.~B}, 43:187--222, 1987.

\bibitem{lopl86}
L.~Lov\'asz and M.~D. Plummer.
\newblock {\em Matching Theory}.
\newblock Number~29 in Annals of Discrete Mathematics. Elsevier Science, 1986.

\bibitem{lumu24}
C.~L. Lucchesi and U.~S.~R. Murty.
\newblock {\em Perfect Matchings: A Theory of Matching Covered Graphs}.
\newblock Springer, 2024.

\bibitem{mccu98}
W.~McCuaig.
\newblock Even {D}icycles.
\newblock {\em J.~Graph Theory}, 35:46--68, 09 2000.

\bibitem{mccu04}
W.~McCuaig.
\newblock P{\'o}lya's {P}ermanent {P}roblem.
\newblock {\em The Electronic J.~of Combin.}, 11, 2004.

\bibitem{mrst97}
W.~McCuaig, N.~Robertson, P.~D. Seymour, and R.~Thomas.
\newblock Permanents, {P}faffian orientations and even directed circuits.
\newblock {\em STOC '97 Proceedings of the twenty-ninth annual ACM symposium on
  Theory of computing}, pages 402--405, 1997.

\bibitem{nadd82}
D.~Naddef.
\newblock Rank of maximum matchings in a graph.
\newblock {\em Mathematical Programming}, 22:52--70, 1982.

\bibitem{rst99}
N.~Robertson, P.~D. Seymour, and R.~Thomas.
\newblock Permanents, {P}faffian orientations and even directed circuits.
\newblock {\em Ann.~of Math.~(2)}, 150:929--975, 1999.

\bibitem{szig02}
Z.~Szigeti.
\newblock Perfect matchings versus odd cuts.
\newblock {\em Combinatorica}, 22:575--589, 2002.

\bibitem{thom06}
R.~Thomas.
\newblock A survey of {P}faffian orientations of graphs.
\newblock 2006.

\bibitem{vali79}
L.~Valiant.
\newblock The complexity of computing the permanent.
\newblock {\em Theoretical Computer Science}, 8:189--201, 1979.

\bibitem{vaya89}
V.~V. Vazirani and M.~Yannakakis.
\newblock Pfaffian orientations, 0--1 permanents, and even cycles in directed
  graphs.
\newblock {\em Discrete Applied Math.}, 25:179--190, 1989.

\end{thebibliography}

\end{document}